\documentclass[a4paper]{article}

\usepackage{amssymb, amsmath, amsthm}
\usepackage{bbm,dsfont}
\numberwithin{equation}{section}
\title{On the 1D Cubic Nonlinear Schr\"odinger Equation in an Almost Critical Space}
\author{Shaoming Guo}
\date{}

\def\R{\mathbb{R}}
\def\N{\mathbb{N}}

\def\Z{\mathbb{Z}}

\def\mb{\mathbb}

\def\mlimits{\displaystyle\sum\limits}
\def\msup{\displaystyle\sup}
\def\mint{\displaystyle\int}
\def\emb{\hookrightarrow}
\def\noin{\noindent}
\def\any{\forall}
\def\del{\Delta}
\def\v{\Vert}

\def\lam{\lambda}
\def\lesim{\lesssim}
\def\begineq{\begin{equation}}
\def\endeq{\end{equation}}

\theoremstyle{plain}
\newtheorem{thm}{Theorem}[section]
\newtheorem{prop}[thm]{Proposition}
\newtheorem{lem}[thm]{Lemma}
\newtheorem{cor}[thm]{Corollary}
\newtheorem{defi}[thm]{Definition}
\newtheorem{claim}[thm]{Claim}
\newtheorem{rem}[thm]{Remark}

\newtheorem*{openproblem*}{Open Problem}

\begin{document}

\maketitle

\begin{abstract}
We consider the Cauchy problem for the one dimensional cubic nonlinear Schr\"odinger equation $iu_t+u_{xx}-|u|^2u=0$. As the first step local well-posedness in the modulation space $M_{2,p}$ ($2\le p<\infty$) is derived (see Theorem \ref{thm1.4}), which covers all the subcritical cases. Afterwards in order to approach the endpoint case, we will prove the almost global well-posedness in some Orlicz type space (see Theorem \ref{thm1.8}), which is a natural generalization of $M_{2,p}$, and is almost critical from the viewpoint of scaling. The new ingredient is an endpoint version of the two dimensional restriction estimate (see Lemma \ref{lem3.6}).
\end{abstract}


\let\thefootnote\relax\footnote{Date: \date{\today}; MSC class: 35Q55}

\section{Introduction and statement of the main result}
In this paper we consider the cubic nonlinear Schr\"odinger equation in dimension one:\\
\begineq\label{main}
\left\{ \begin{array}{ll}								
iu_t+u_{xx}-|u|^2u=0 & (t,x)\in \mb{R}^+\times \mb{R}\\
u(0,x)=u_0(x) & x\in \mb{R}
\end{array} \right.
\endeq

\noin which has been proven to be globally well-posed in the Sobolev spaces $H^s$ ($\any s\ge 0$) by Tsutsumi in \cite{TS}, see also Cazenave and Weissler \cite{CW} for more general nonlinear terms and initial data.\\

This equation has a rich structure, and we will just mention three that are relevant to us:\\
\indent {\bf Scaling:} If $u(t,x)$ is a solution for \eqref{main}, with initial data $u_0(x)$, then 
$$v(t,x):= \lambda u(\lambda^2 t, \lam x)$$
is also a solution for \eqref{main} with the corresponding initial data $v_0(x):=\lambda u_0(\lambda x)$. For the Sobolev spaces, the scaling invariant index is $s=-1/2$, i.e. 
\begineq
\v u_0(x)\v_{\dot{H}^{-1/2}}=\v v_0(x)\v_{\dot{H}^{-1/2}},
\endeq
which suggests that we might expect well-posedness result for all initial data in $H^s$ with $s>-1/2$.\\
\indent {\bf Galilean invariance:} Under the above notations, we can check easily that 
\begineq
u_c(t,x):=e^{-i(c^2 t-cx)}u(t, x-2ct)
\endeq
is also a solution for \eqref{main} with initial data $u_c(0,x):=e^{icx}u_0(x)$.\\
\indent By applying the Galilean invariance (which implies the existence of any possible speed of propagation), Kenig, Ponce and Vega in \cite{KPV} derived the ill-posedness of \eqref{main} with initial data in $H^s$, for $ s\in [-1/2, 0)$. \\
\indent {\bf Complete integrability:} This equation is completely integrable, so there are infinitely many conservation laws. By the inverse scattering method, the asymptotics of the solution with Schwartz initial data can be obtained. We refer to Deift and Zhou \cite{DZ} for the precise statement of the results and the references therein.\\

As we can see that there is 1/2 derivative gap left between the well-posedness result suggested by the scaling ($\dot{H}^s$, with $s=-1/2$), and the one that has been proved (well-posedness in $H^s$ with $s\ge 0$, and ill-posedness in $H^s$ with $s<0$), in order to understand the behavior of the solution under $L^2$, a lot of work has been done. In \cite{CCT2}, \cite{KT2}, \cite{KT3}, a priori bound for solutions with the initial data in $H^s$ for some negative $s$ is obtained. In terms of well-posedness, in Vargas and Vega \cite{VV}, initial data with infinite $L^2$ norm is considered, and the local well-posedness is achieved for initial data satisfying 
\begineq
\v e^{it \partial_x^2} u_0\v_{L_t^3L_x^6(I\times \mb{R})}<+\infty.
\endeq
Examples for this kind of initial data are also given there, which read as:
$\hat{u}_0(\xi)\in C^1$, for some $\beta>1/6$, and for $j=0, 1$, 
\begineq\label{E:1.2}
|\frac{d^j}{d\xi^j}\hat{u}_0(\xi)|\le \frac{C}{(1+|\xi|)^{\beta+j}}.
\endeq
However, from the viewpoint of scaling, local well-posedness for all initial data satisfying \eqref{E:1.2} with any $\beta>0$ can be expected.\\

By generalizing the Bourgain space method and working directly in the frequency space, Gr\"{u}nrock in \cite{Gr} derived the well-posedness result for more general initial data. Remarkably for the local well-posedness issue, the 1/2 gap of derivative stated above was covered for the first time. Let us state briefly the result there.
\begin{thm}(\cite{Gr}) For all $r\in (1, \infty)$, \eqref{main} is locally well-posed for initial data in $\hat{L}^r(\mb{R})$, with the norm given by
\begineq
\v u_0\v_{\hat{L}^r(\mb{R})}:=\v \hat{u}_0\v_{L^{r'}(\mb{R})},
\endeq
where $r'$ is the dual number of $r$, i.e. $1/r+1/r'=1.$
\end{thm}

\begin{rem} The gap left in the result of Vargas and Vega in \cite{VV} is covered by this theorem, as we can easily see that for any initial data satisfying 
$$|\hat{u}_0(\xi)|\le \frac{C}{(1+|\xi|)^{\beta}}, \any \beta>0,$$
there exists $r>1$, such that $u_0\in \hat{L}^r(\mb{R})$.
\end{rem}

\begin{rem} The scaling of $\hat{L}^r$. After some standard calculation we can see that the space $\hat{L}^r$ has the same scaling as $\dot{H}^s$, $s=1/2-1/r$. Notice that when $r\to 1^+$, we have $s\to -\frac{1}{2}^+$, and that $\dot{H}^{-1/2}$ and $\hat{L}^1$ are both the scaling invariant spaces. We say in this sense that the scaling invariant index is almost reached.
\end{rem}

Now we want to generalize the result of Gr\"{u}nrock, and work in a relatively larger space for initial data. To be precise, we will prove:
\begin{thm}\label{thm1.4}
 The equation \eqref{main} is locally well-posed for initial data in the modulation spaces $M_{2,p}$ with $2\le p<\infty$.
\end{thm}

Let us recall how the norm of $M_{2,p}$ is defined. For $k\in \mathbb{Z}$, denote $\mathbbm{1}_{[k,k+1]}$ as the characteristic function of the interval $[k,k+1]$, define the frequency projection operator $P_k$ by 
\begineq
P_k(u)=\mathcal{F}^{-1}\mathbbm{1}_{[k,k+1]}(\xi)\mathcal{F}u,
\endeq
then 
\begineq\label{2608ee1.8}
\v u_0\v_{M_{2,p}}:=(\mlimits_{k\in \mb{Z}} \v P_k(u_0)\v_{L^2}^p)^{1/p}.
\endeq
\begin{rem}
We have the inclusion of the function spaces 
\begineq\label{E:1.3}
\hat{L}^{p'}\emb M_{2,p},
\endeq
which follows from
\begin{eqnarray*}
\v u\v_{M_{2,p}} 	&=&(\sum_{k\in \mb{Z}} \v \hat{u}_0(\xi)\mathbbm{1}_{[k,k+1]}(\xi)\v_{L^2}^p)^{1/p}\\
			&\le& C(\sum_{k\in \mb{Z}} \v \hat{u}_0(\xi)\mathbbm{1}_{[k,k+1]}(\xi)\v_{L^p}^p)^{1/p}\\
			&\le& C(\int_{\mb{R}}|\hat{u}_0(\xi)|^p d\xi)^{1/p}\le C\v u_0\v_{\hat{L}^{p'}}.\\	
\end{eqnarray*}
\end{rem}

\begin{rem} The endpoint case $M_{2,\infty}$. A typical example for initial data in this space is the Dirac function $\delta_0(x)$. There will be a  serious problem in the uniqueness of the solution when solving \eqref{main} with initial data $u_0(x)=\delta_0(x)$. We refer to Theorem 1.5 in Kenig, Ponce and Vega \cite{KPV} for further details.
\end{rem}

\begin{rem} As here we only consider the local well-posedness issue, it does not matter if the equation is focusing or defocusing. While for the ill-posedness in the defocusing case, due to the lack of soliton solutions, the argument is different from the one given by Kenig, Ponce and Vega, see Christ, Colliander and Tao \cite{CCT}, where the pseudo-conformal transformation is applied.\\
\end{rem}

Similar to \eqref{2608ee1.8}, one can define $M_{q, p}$ for all $p, q\ge 1$. These spaces are called modulation spaces. Modulation spaces were introduced by Feichtinger in \cite{Feich}. Recently they have been used intensively in the study of nonlinear dispersive equations. This was initiated by Wang and his collaborators, see \cite{WH1}, \cite{WH2} and \cite{WZG}. The results therein include well-posedness results for various dispersive equations with initial data in $M_{q, 1}$ for certain $q\in (1, \infty)$. Benyi and Okoudjou \cite{BO} generalised (part of) the above results to $M_{q, 1}$ for all $q\ge 1$. We refer to the survey paper by Ruzhansky, Sugimoto and Wang \cite{RSW} and the book by Wang et al. \cite{WHHG} for more results in the same spirit. 

To the author's knowledge, the result in the above Theorem \ref{thm1.4} is the first time that a nonlinear Schr\"odinger equation is proved to be locally well-posed for initial data in a modulation space $M_{2, p}$ for certain $p>1$. Indeed, our result covers the whole possible range of $p$.

Another closely related study is on the Schr\"odinger equations with rough potentials in modulations spaces. We refer to \cite{CGNR}, \cite{CN}, \cite{CNR}. The main techniques used there, for example the so-called Gabor decomposition, are different from \cite{BO}, \cite{WH1}, \cite{WH2}, \cite{WZG} and the present paper. For more results on dispersive equations in modulation spaces, we refer to \cite{BGOR}, \cite{CN0}, \cite{KKI1} and \cite{KKI2}.\\

We proceed with presenting our results. With in mind two examples $\hat{u}_0(\xi)=\frac{1}{(1+|\xi|)^{\beta}},\any \beta>0$, for which we have the local well-posedness, and $\hat{u}_0(\xi)=1$, for which we do not, a natual question is to see how close we can approach the endpoint case. Now we state the second theorem, from which the almost global well-posedness follows:

\begin{thm}\label{thm1.8}
The equation \eqref{main} is well-posed on the time interval $[0,1]$ for small initial data in $l^{\Phi}L^2$, where the norm is given by
\begineq\label{orlicz}
\v u_0\v_{l^{\Phi}L^2}:=\inf\{\lambda>0: \sum_{n\in \mb{Z}} \Phi(\frac{\v P_n(u_0)\v_{L^2}}{\lambda})\le 1\},
\endeq
with 
\begin{displaymath}
\Phi (x)=\left\{ \begin{array}{ll}								
e^{-(1/x)^{1/\gamma}+ C_{\gamma}x}& {x>0}\\
0 & {x=0}
\end{array} \right.
\end{displaymath}
for any $\gamma>2$, and a large enough constant $C_{\gamma}$. See the remark below for the choice of the constant $C_{\gamma}$.
\end{thm}

\begin{rem} 
One typical example for initial value in the function space $l^{\Phi}L^2$ is 
\begineq\label{initial}
|\hat{u}_0(\xi)|\sim \frac{1}{\ln^{\gamma}(2+|\xi|)},
\endeq
which decays slower than $\frac{1}{(1+|\xi|)^{\beta}}$ for any $\beta$ when $\xi\to \infty$.
\end{rem}

\begin{rem}The norm we use in \eqref{orlicz} is the discrete Orlicz norm with respect to the function $\Phi$, which is a natural generalization of the $M_{2,p}$ space.
\end{rem}

\begin{rem}
Clearly the main decaying part for $\Phi$ when $x\to 0^+$ is $e^{-(1/x)^{1/\gamma}}$, here we use $e^{-(1/x)^{1/\gamma}+C_{\gamma} x}$ instead of $e^{-(1/x)^{1/\gamma}}$ just because we need $\Phi$ to be convex in the definition of the Orlicz space (see Definition \ref{defi5.1} and \ref{defi5.2} below), which can be guaranteed by choosing $C_{\gamma}$ large enough.
\end{rem}

\begin{rem}
Similar to the embedding \eqref{E:1.3}, in the case of the Orlicz space, we have
\begineq\label{E:1.5}
\hat{L}^{\Phi}\emb l^{\Phi}L^2,
\endeq
where the norm for the function space $\hat{L}^{\Phi}$,which is a natural generalization of $\hat{L}^r$, is given by \\
\begineq
\v u_0\v_{\hat{L}^{\Phi}}:= \inf\{\lam>0, \int_{\mb{R}}\Phi(\frac{|\hat{u}_0(\xi)|}{\lam})d\xi \le 1\}.
\endeq
The proof for \eqref{E:1.5} is similar to that of \eqref{E:1.3}, hence we leave it out.
\end{rem}

\begin{cor}\label{cor1.13}
The equation \eqref{main} is almost globally well-posed for small initial data in $\hat{L}^{\Phi}\cap \hat{L}^1$, with the time span $T$ given by\\
$$e^{(\v u_0\v_{\hat{L}^{\Phi}}+\v u_0\v_{\hat{L}^1})^{-1/\gamma}},$$
and the solution u lies in the space $X_{\Phi}$, which is
$$\v u\v_{X_{\Phi}}:=\v\v P_n(u)\v_{U_{\del, [0,T]}^2}\v_{l_n^{\Phi}}<+\infty,$$
see Section 2 for the definition of the function space $U_{\del}^2$. Then by the embedding relation $U_{\del}^2\emb L^{\infty}(\mb{R}, L^2)$ in Proposition \ref{prop2.2}, we see easily that 
$$\v\v\v P_n(u)\v_{L_x^2}\v_{l_n^{\Phi}}\v_{L_t^{\infty}}\lesim \v\v P_n(u)\v_{L_t^{\infty}L_x^2}\v_{l_n^{\Phi}} <+\infty,$$
i.e. the $l^{\Phi}L^2$ regularity of the initial data is persisted.\\
\end{cor}



{\bf Notations:} For $f, g$ being two non-negative functions, by ``$f\lesim g$'' we mean that there exists a constant $C$ such that $f(x)\le C g(x), \any x\in\bf{R}$.

For $a, b\in \bf{R}$, by $a\sim b$ we mean that $|a-b|$ is less than some universal constant, say 4. Similarly by $a\ll b$ we mean $a-b\le 4$, and $a\gg b$ means $a-b\ge 4$.

We also need another concept of ``being comparable''. Take a universal constant $C>1$, for $a,b\in \mb{R^{+}}$, by $a \approx b$ we mean that $1/C<a/b<C$.

The space-time norm will be used later in a large scale. Hence for the sake of simplicity, we want to make the convention that whenever the norm like $L^pL^q$ appears, we always mean  $L_t^pL_x^q$, if not stated otherwise.

For an interval $I\subset \mb{R}$, by $P_{I}$ we mean the usual frequency projection operator on $I$, i.e. $P_{I}(u):=\mathcal{F}^{-1}\mathbbm{1}_{I}\mathcal{F}(u).$

For $n\in \mathbb{Z}$, we have defined $P_n$ when introducing the norm of the modulation spaces. Sometimes we will also use the notation $u_n$ instead of $P_n(u)$ for simplicity.\\

{\bf Organization of paper:} In the second section, we will introduce the function spaces that we will work in, which first appeared in Koch and Tataru \cite{KT1}, and afterwards were characterised in detail in Koch and Tataru \cite{KT2}, and Hadac, Herr and Koch \cite{HHK}.

In the third section, we state the known linear and bilinear estimates, and prove a crucial lemma (Lemma \ref{lem3.6}), which will be useful in proving the trilinear estimates.

In the fourth section, the proof of the trilinear estimate in modulation space $M_{2,p}$ will be given, then Theorem \ref{thm1.4} follows directly from the standard contraction mapping principle.

In the fifth section, we will give a brief introduction to the Orlicz spaces, and collect some technical lemmas which we will need in the proof of trilinear estimate in the Orlicz spaces.

In the last section, which is totally parallel to Section 4, we will give the proof of the trilinear estimate in the above defined Orlicz type space, and again by the standard contraction mapping principle, Theorem \ref{thm1.8} follows.\\

{\bf Acknowledgements:} This work is done under the supervision of Prof. Herbert Koch. The author would like to thank him for his patient guidance and for sharing many helpful thoughts. The author also thanks Prof. Sebastian Herr for carefully reading this paper and giving a lot of valuable suggestions. The author would also like to thank a anonymous referee for pointing out several references, and other valuable comments.

\section{The $U^p$ and $V^p$ spaces}

Most of the materials in this section can be found in \cite{HHK}. For the sake of completeness, we still include them here.

\begin{defi}\label{defi2.1} (\cite{HHK}) Let $\mathcal{Z}$ be the set of finite partitions $-\infty =t_0<t_1<\dots<t_K=\infty$. Let $1\le p < \infty$, for $\{t_k\}_{k=0}^K \in \mathcal{Z}$ and $\{\phi_k\}_{k=0}^{K-1} \subset L^2$ with $\sum_{k=0}^{K-1} \v \phi_k\v_{L^2}^p=1$ and $\phi_0=0$ we call the function $a: \mb{R}\to L^2$ given by
\begineq\label{atom}
a=\sum_{k=1}^K \chi_{[t_{k-1}, t_k)}\phi_{k-1}
\endeq
a $U^p$-atom. Furthermore, we define the atomic space 
\begineq
U^p:= \{u=\sum_{j=1}^{\infty} \lam_j a_j \big\vert \text{with } a_j \text{ some } U^p\text{-atom}, \lam_j \in \mb{C} \text{ such that} \sum_{j=1}^{\infty}|\lam_j|<\infty\}
\endeq
with norm 
\begineq\label{unorm}
\v u\v_{U^p}:=\inf\{\sum_{j=1}^{\infty}|\lam_j| \big\vert u=\sum_{j=1}^{\infty}\lam_j a_j,  \lam_j\in \mb{C}, a_j\text{ is } U^p\text{-atom} \}.
\endeq
\end{defi}

\begin{prop}\label{prop2.2}(\cite{HHK})
Let $1\le p<q<\infty$,\\
(i) $U^p$ is Banach space;\\
(ii) The embedding $U^p\subset U^q\subset L^{\infty}(\mb{R}, L^2)$ is continuous;\\
(iii) For $u\in U^p$ it holds $\lim_{t\to t_0^+} \v u(t)-u(t_0)\v_{L^2}=0$, i.e. every $u\in U^p$ is right continuous.
\end{prop}

\begin{defi}\label{defi2.3} (\cite{HHK},\cite{Wi})
Let $1\le p<\infty$, the space $V^p$ is defined as the normed space of all functions $v: \mb{R}\to L^2$ such that $\lim_{t\to \pm \infty}v(t)$ exist and for which the norm 
\begineq
\v u\v_{V^p}:=\sup_{\{t_k\}_{k=0}^K \in \mathcal{Z}}(\sum_{k=1}^K \v v(t_k)-v(t_{k-1})\v_{L^2}^p)^{1/p}
\endeq
is finite, where we use the convention $v(-\infty)=\lim_{t\to -\infty}v(t)$ and $v(\infty)=0$. Let $V_-^p$ denote the closed subspace of all $v\in V^p$ with $\lim_{t\to -\infty}=0$.
\end{defi}

\begin{prop}\label{prop2.4}(\cite{HHK})
 Let $1\le p<q<\infty$,\\
(i) we define the closed subspace $V_{rc}^p(V_{-,rc}^p)$ of all right continuous $V^p$ functions ($V_-^p$ functions). The spaces $V^p, V_{rc}^p, V_-^p$ and $V_{-,rc}^p$ are Banach spaces;\\
(ii) the embedding $U^p\subset V_{-,rc}^p$ is continuous;\\
(iii) the embeddings $V^p\subset V^q$, $V_-^p\subset V_-^q$ and $V_{-,rc}^p\subset V_{-,rc}^q$ are continuous.
\end{prop}

\begin{prop}\label{prop2.5}(Interpolation) For $1\le p<q<\infty$, there exists a positive constant $\epsilon(p,q)>0$, s.t. $\any u\in V^p$, $\any M\in \mb{Z}^+$, there exists $u_1, u_2,$ s.t. $u=u_1+u_2$, and 
\begineq
\frac{1}{M}\v u_1\v_{U^p}+e^{\epsilon M}\v u_2\v_{U^q}\lesim \v u\v_{V^p}.
\endeq
Moreover we could see that the embedding $V_{-, rc}^p\subset U^q$ is continuous.
\end{prop}

\noin {\bf Proof of Proposition \ref{prop2.5}:} See the proof of Proposition 2.5 and Corollary 2.6 in\cite{HHK}.$\square$

\begin{prop} \label{prop2.6}(orthogonality in $U^2$ and $V^2$)\\ 
1) Take an interval $I:=[m, n)\subset \mb{R}$ with $m,n\in \mb{Z}$, then for $u\in V^2$ the following orthogonality holds:
\begineq\label{vortho}
\v P_{I}u\v_{V^2}\le (\mlimits_{i\in I, i\in \mb{Z}} \v P_i u\v_{V^2}^2)^{1/2}.
\endeq
2) Similarly for $u\in U^2$, the following orthogonality holds:
\begineq\label{uortho}
\v P_{I}u\v_{U^2}\ge (\mlimits_{i\in I, i\in \mb{Z}} \v P_i u\v_{U^2}^2)^{1/2}.
\endeq
\end{prop}

\noindent {\bf Proof of Proposition \ref{prop2.6}:} For the estimate \eqref{vortho},
\begin{eqnarray*}
\v P_I u\v_{V^2}^2 &=& \msup_{\{t_k\}_{k=0}^K \in \mathcal{Z}}\mlimits_{k=1}^K \v (P_I u)(t_k)-(P_I u)(t_{k-1})\v_{L^2}^2\\
			&\le& \msup_{\{t_k\}_{k=0}^K \in \mathcal{Z}}\mlimits_{k=1}^K \mlimits_{i\in I, i\in \mb{Z}} \v (P_i u)(t_k)-(P_i u)(t_{k-1})\v_{L^2}^2\\
			&\le& \mlimits_{i\in I, i\in \mb{Z}} \msup_{\{t_k\}_{k=0}^K \in \mathcal{Z}}\mlimits_{k=1}^K \v (P_i u)(t_k)-(P_i u)(t_{k-1})\v_{L^2}^2\\
			&\le & \mlimits_{i\in I, i\in \mb{Z}} \v P_i u\v_{V^2}^2
\end{eqnarray*}
where the first step is by the definition of $V^2$ space, the second is by the orthogonality in $L^2$, and the last is again the definition of $V^2$ space.\\
\indent For the estimate \eqref{uortho}, take $u\in U^2$, notice that $P_I$ is a projection operator, i.e. $P_I P_I=P_I$, so it's enough to consider $u$ with $\text{supp} \mathcal{F}(u)\subset I$. By the definition of $U^2$ space, $\any \epsilon>0$, there exists a representation
\begineq
u=\mlimits_{j=1}^{\infty}\lam_j a_j,
\endeq
with $\lam_j\in \mb{C}$, $a_j$ some $U^2\text{-atom},$ such that 
\begineq
\mlimits_{j=1}^{\infty}|\lam_j|\le \v u\v_{U^2}+\epsilon.
\endeq
Then 
$$(\mlimits_{i\in I, i\in \mb{Z}} \v P_i u\v_{U^2}^2)^{1/2}\lesim \mlimits_{j=1}^{\infty}|\lam_j| (\mlimits_{i\in I, i\in \mb{Z}} \v P_i a_j\v_{U^2}^2)^{1/2}\lesim \mlimits_{j=1}^{\infty}|\lam_j|\lesim \v u\v_{U^2}+\epsilon,$$
as $\epsilon$ can be chosen arbitrarily small, we come to the conclusion.$\square$

\begin{prop}\label{prop2.7}(\cite{KT2},\cite{HHK}) Let $1<p<\infty$, we have 
\begineq
(U^p)^*=V^{p'}
\endeq
in the sense that 
\begineq
\begin{split}
T: V^{p'}\to (U^p)^*,\\
T(v):=B(\cdot, v),
\end{split}
\endeq
is an isometric isomorphism, where the bilinear form $B: U^p\times V^{p'}\to \mb{R}$ is defined in the following way: first for a partition $\mathfrak{t}:=\{t_k\}_{k=0}^K\in \mathcal{Z}$ we define
\begineq
B_{\mathfrak{t}}(u,v):=\sum_{k=1}^K \langle u(t_{k-1}), v(t_k)-v(t_{k-1})\rangle.
\endeq
Here $\langle \text{ } ,\text{ } \rangle$ denotes the $L^2$ inner product on $\R$. Then for any $u\in U^p, v\in V^{p'}$, there exists a unique number $B(u,v)$ such that $\any \epsilon>0$ there exists a partition $\mathfrak{t}\in \mathcal{Z}$ such that $\any \mathfrak{t}'\supset \mathfrak{t}$ it holds that 
$$|B_{\mathfrak{t}'}(u,v)-B(u,v)|<\epsilon.$$
\end{prop}

\begin{defi}\label{defi2.8}
We define \\
(i) $U_{\del}^p=e^{\cdot i\del}U^p$ with norm $\v u\v_{U_{\del}^p}:=\v e^{-it\del}u\v_{U^p},$\\
(ii) $V_{\del}^p=e^{\cdot i\del}V^p$, with norm $\v u\v_{V_{\del}^p}:=\v e^{-it\del}u\v_{V^p},$ and similarly the closed subspaces $V_{rc, \del}^p, V_{-, \del}^p, V_{-, rc, \del}^p$.
\end{defi}

\begin{rem}
In the following setting, we will always consider the right continuous functions, so we will just write $V_{\del}^p$ instead of $V_{-, rc, \del}^p$ for simplicity.
\end{rem}

\begin{rem}\label{rem2.10}
Later in the proof of the trilinear estimates, we need this duality relation in Proposition \ref{prop2.7} because we will estimate norms like 
\begineq\label{duality}
\v \int_0^t e^{i(t-s)\del}F(x,s)ds\v_{U_{\del}^2}.
\endeq
Let us show how to simplify the above expression:
\begin{eqnarray*}
\eqref{duality}	&=&\v \mint_0^t e^{-is\del}F(x,s)ds\v_{U^2}\\
			&=&\sup_{v\in V^2, \v v\v_{V^2}\le 1} B(\mint_0^t e^{-is\del}F(x,s)ds, v(x,t))\\
			&=&\sup_{v\in V^2, \v v\v_{V^2}\le 1} \mint_{\bf{R}}\langle e^{-it\del}F(x,t),v(x,t)\rangle dt\\
			&=&\sup_{v\in V^2, \v v\v_{V^2}\le 1} \mint_{\bf{R}}\langle F(x,t),e^{it\del}v(x,t)\rangle dt\\
			&=&\sup_{v\in V^2, \v v\v_{V_{\del}^2}\le 1} \mint_{\bf{R}}\langle F(x,t),v(x,t)\rangle dt,
\end{eqnarray*}
where the first equality is by the definition of $U_{\del}^2$, the second is by the duality relation given in Proposition \ref{prop2.7}, the third is basically integration by part with respect to the time variable $t$, and the last is by the definition of $V_{\del}^2$.
\end{rem}

\begin{prop}\label{prop2.9}(\cite{KT2})
Define the homogeneous Besov type Bourgain space $\dot{X}^{s,1/2,1}$ and $\dot{X}^{s,1/2,\infty}$ with the norm 
\begineq
\v u\v_{\dot{X}^{s,1/2,1}}:= \sum_{\mu \text{ dyadic}}(\int_{|\tau-\xi^2|\approx \mu}|\hat{u}(\tau, \xi)|^2 |\xi|^{2s}|\tau-\xi^2|d\xi d\tau)^{1/2}
\endeq
\begineq
\v u\v_{\dot{X}^{s,1/2,\infty}}:= \sup_{\mu \text{ dyadic}}(\int_{|\tau-\xi^2|\approx \mu}|\hat{u}(\tau, \xi)|^2 |\xi|^{2s}|\tau-\xi^2|d\xi d\tau)^{1/2},
\endeq
then the following embedding holds
\begineq
\dot{X}^{0,1/2,1}\subset U_{\del}^2\subset V_{\del}^2\subset \dot{X}^{0,1/2,\infty}
\endeq
\end{prop}

\begin{rem} 
Later we will apply this proposition in such a way that, if the modulation of a  function $u$ is no less than $\mu$, which is some dyadic number, then
\begineq
\v u\v_{L^2L^2}\lesim \mu^{-1/2}\v u\v_{V_{\del}^2}
\endeq
which can be easily seen from the above embedding.
\end{rem}

\section{Linear and Bilinear Estimates}

\begin{prop}\label{prop3.1} (Strichartz estimate, \cite{Str})
Let $p, q$ be indices satisfying 
\begineq
\frac{2}{p}+\frac{1}{q}=\frac{1}{2}, 4\le p\le \infty
\endeq
then the solution of the homogeneous equation
\begineq\label{free}
iu_t-\del u=0
\endeq
with initial data $u(0,x)=u_0(x)$ satisfies the Strichartz estimate:\\
$$\v u\v_{L^p_t L^q_x}\lesim \v u_0\v_{L^2}.$$
\end{prop}


\begin{cor}\label{cor3.2}
With $p,q$ defined above, we have the following estimate concerning the norm which we will use later:
$$\v v\v_{L_t^pL_x^q}\lesim \v v\v_{U_{\del}^p}.$$ 
\end{cor}
\noindent {\bf Proof of Corollary \ref{cor3.2}:} The proof is quite straightforward, as it suffices to check the estimate only for atoms, which then is just the Strichartz estimate in Proposition \ref{prop3.1}.$\square$

\begin{prop}\label{prop3.3}(Bilinear estimate; \cite{Gr}, \cite{KT3})
Let $\lam>0$, assume that $u, v$ are solutions to the linear equation \eqref{free} with the corresponding initial data $u_0$ and $v_0$, then 
\begineq
\v P_{>\lam}(u \bar{v})\v_{L^2L^2}\lesim \lam^{-1/2}\v u_0\v_{L^2}\v v_0\v_{L^2}
\endeq
\end{prop}

\noin {\bf Proof of Proposition \ref{prop3.3}:} This estimate is also well-known. As the calculation is not that involved, we will do it here. By definition of $u,v$, we have 
$$u=\mathcal{F}_x^{-1} e^{-it\xi^2}\mathcal{F}_x u_0, v=\mathcal{F}_x^{-1} e^{-it\xi^2}\mathcal{F}_x v_0.$$
then 
\begin{eqnarray*}
\mathcal{F}_x (u \bar{v})&=& (e^{-it\xi^2}\hat{u}_0(\xi))*(e^{it\xi^2}\bar{\hat{v}}_0(\xi))\\
				&=& \int_{\mb{R}} e^{-it\xi^2+2it\xi\xi_1}\hat{u}_0(\xi-\xi_1)\bar{\hat{v}}_0(\xi_1)d\xi_1
\end{eqnarray*}
\begin{eqnarray*}
\mathcal{F}_{x,t}(u\bar{v})&=&\int_{\mb{R}}\delta_{\xi^2-2\xi\xi_1+\tau}\hat{u}_0(\xi-\xi_1)\bar{\hat{v}}_0(\xi_1)d\xi_1\\
		&=&\frac{1}{2|\xi|}\hat{u}_0(\frac{\xi^2-\tau}{2\xi})\hat{v}_0(\frac{\xi^2+\tau}{2\xi})
\end{eqnarray*}
now calculate the $L^2_{x,t}$ norm:
\begin{eqnarray*}
\v P_{\lam}(u\bar{v})\v_{L^2}^2 &=&\int_{\mb{R}^2} \mathbbm{1}_{\{|\xi|>\lam\}}\frac{1}{4\xi^2} |\hat{u}_0(\frac{\xi^2-\tau}{2\xi})|^2 |\hat{v}_0(\frac{\xi^2+\tau}{2\xi})|^2 d\xi d\tau\\
		&=& \int_{\mb{R}^2} \mathbbm{1}_{\{|\xi_1+\xi_2|>\lam\}}\frac{1}{|\xi_1+\xi_2|}|\hat{u}_0(\xi_1)|^2|\hat{v}_0(\xi_2)|^2d\xi_1 d\xi_2
\end{eqnarray*}
then the conclusion follows.$\square$\\

Again by testing atoms in $U_{\Delta}^2$ space, the above proposition implies
\begin{cor}\label{cor3.4}
The following estimate holds:
\begineq\label{E3.4}
\v P_{>\lam}(u \bar{v})\v_{L_{x,t}^2}\lesim \lam^{-1/2}\v u\v_{U_{\del}^2}\v v\v_{U_{\del}^2}
\endeq
\end{cor}

Combining the above bilinear estimate with the interpolation of $U^2$ and $V^2$ spaces in Proposition \ref{prop2.5} we get the following variant of the bilinear estimate:
\begin{cor}\label{cor3.5}
For $m,n\in\mathbb{Z}$, with $m\gg n$, and $u_m, u_n$ the frequency localization of $u$ near $m$ and $n$, the following estimate holds on the time interval $[0, 1]$:
\begineq\label{interbilinear}
\v u_m u_n\v_{L^2L^2}\lesim \frac{\ln^2(m-n)}{(m-n)^{1/2}}\v u_m\v_{V_{\del}^2}\v u_n\v_{V_{\del}^2}
\endeq
Moreover for the proof of the large data local well-posedness, we will work in a small time interval $[0,T]$, then the following estimate becomes necessary:
\begineq\label{interbilinear2}
\v u_m u_n\v_{L^2L^2}\lesim T^{\theta/4} \frac{\ln^{2(1-\theta)}(m-n)}{(m-n)^{(1-\theta)/2}}\v u_m\v_{V_{\del}^2}\v u_n\v_{V_{\del}^2}
\endeq
where $\theta$ is any given real number in $(0,1)$.
\end{cor}

\begin{rem}
Here we do interpolation for two terms, i.e. replace both $U_{\del}^2$ norms on the right hand side of \eqref{E3.4} by $V_{\del}^2$. In fact it suffices to do only one, but this will not give any improvement to the final result (for example, this does not improve the restriction $\gamma>2$ in Theorem \ref{thm1.8}), hence we choose to do both, which will also simplify the forthcoming calculation.
\end{rem}

\noin{\bf Proof of Corollary \ref{cor3.5}:} The first estimate will be obtained by interpolating \eqref{E3.4} with the following H\"older's inequality 
\begineq
\|u_m u_n\|_{L^2 L^2}\lesim \|u_m\|_{L^4L^4}\|u_n\|_{L^4L^4}\lesim \|u_m\|_{L^8L^4}\|u_n\|_{L^8L^4}\lesim \|u_m\|_{U_{\del}^8}\|u_n\|_{U_{\del}^8}.
\endeq
To be precise, we take $M\in \mb{N}$ to be chosen later. By Proposition \ref{prop2.5}, there exists $u_{m,1}, u_{m,2},u_{n,1}, u_{n,2}$ such that
$$\frac{1}{M}\v u_{m,1}\v_{U_{\del}^2}+e^{\epsilon M}\v u_{m,2}\v_{U_{\del}^8}\lesim \v u_m\v_{V_{\del}^2},$$
$$\frac{1}{M}\v u_{n,1}\v_{U_{\del}^2}+e^{\epsilon M}\v u_{n,2}\v_{U_{\del}^8}\lesim \v u_n\v_{V_{\del}^2},$$
then 
\begin{eqnarray*}
\v u_m u_n\v_{L^2L^2}&=&\v (u_{m,1}+u_{m,2})(u_{n,1}+u_{n,2})\v_{L^2L^2}\\
				&\lesim& (\frac{M^2}{(m-n)^{1/2}}+M e^{-\epsilon M}+e^{-2\epsilon M})\v u_m\v_{V_{\del}^2}\v u_n\v_{V_{\del}^2}
\end{eqnarray*}
by choosing $M\sim \ln(m-n)$, we could get the bound on the right hand side of $\eqref{interbilinear}$.\\
\indent The estimate $\eqref{interbilinear2}$ just follows from $\eqref{interbilinear}$ and H\"older's inequality:
\begin{eqnarray*}
\v u_m u_n\v_{L^2L^2}&=&\v u_m u_n\v_{L^2L^2}^{\theta}\v u_m u_n\v_{L^2L^2}^{1-\theta}\\
			&\lesim& \v u_m\v_{L^4L^4}^{\theta}\v u_n\v_{L^4L^4}^{\theta}\v u_m u_n\v_{L^2L^2}^{1-\theta}\\
			&\lesim& T^{\theta/4}\v u_m\v_{L^8L^4}^{\theta}\v u_n\v_{L^8L^4}^{\theta}\v u_m u_n\v_{L^2L^2}^{1-\theta}
\end{eqnarray*}
After applying Strichartz estimate and $\eqref{interbilinear}$ we will see the desired estimate.$\square$\\

Now we will state and prove the lemma which will be crucial in the proof of the trilinear estimates:

\begin{lem}\label{lem3.6} 
Take an interval $I= [a, b] \subset \mb{R}$, with $a, b\in \mathbb{Z}$, denote $|I|$ as the length of this interval, define $P_I$ as usual, then \\
1) let $u$ be the solution of the homogeneous equation \eqref{free} on the time interval $[0,1]$ with initial data $u_0$, the the following estimate holds:
\begineq\label{EE3.7}
\v P_I u\v_{L^4([0,1]\times \mb{R})} \lesim \sqrt{\ln|I|} (\mlimits_{k\in I, k\in \mb{Z}} \v u_{0,k}\v_{L^2(\mb{R})}^4)^{1/4}.
\endeq

\noin 2) For any $\beta>1$, similar to the above estimate for the free solution, we have:
\begineq\label{EE3.8}
\v P_I u\v_{L^4([0,1]\times \mb{R})} \lesim (\ln|I|)^{\beta+1/2} (\mlimits_{k\in I, k\in \mb{Z}} \v u_{k}\v_{V_{\del}^2}^4)^{1/4}.
\endeq

\noin 3) In order to gain a time factor in the trilinear estimate, we also need the following variant of the above estimates: $\any \theta\in (0,1)$, $\beta>0$, 
\begineq\label{EE3.9}
\v P_I u\v_{L^4([0,T]\times \mb{R})}^2\lesim 
(T^{\frac{1-\theta}{4}}|I|^{2\beta+\frac{1-\theta}{2}}+T^{1/4})(\mlimits_{k\in I, k\in \mb{Z}} \v u_k\v_{V_{\del}^2}^4)^{1/2}.
\endeq
\end{lem}

\begin{rem}
The estimates above are the endpoint version for the restriction theorem in two dimensions. Define
$$(Tf)(\xi_1, \xi_2):= \int_{\mb{R}} e^{ix\xi_1 +ix^2 \xi_2}\psi(x, \xi_1, \xi_2)f(x)dx$$\\
mapping $f: \mb{R}\to \mb{R}$ to $Tf: \mb{R}^2\to \mb{R}$, where $\psi$ is a cut-off function with compact support in both variables. The restriction theorem then says that 
$$\v Tf\v_{L^q(\mb{R}^2)}\lesim \v f\v_{L^p(\mb{R})},$$
where $q=3p'$, and $1\le p<4$.
Notice that $p=4$ corresponds to $q=4$, i.e. what we derive here is the endpoint version of the two dimensional restriction theorem with a logarithm loss.
\end{rem}


\noin {\bf Proof of Lemma \ref{lem3.6}:} For the estimate in \eqref{EE3.7},
\begin{align*}
& \v P_I u \v_{L^4([0,1]\times \mb{R})}^2 = \v (P_I u)^2 \v_{L^2([0,1]\times \mb{R})}\\
&=\v \mlimits_{m, n\in I, m, n\in \mb{Z}} u_m u_n\v_{L^2([0,1]\times \mb{R})} \\
&=\v \mlimits_{l\in \mb{N}}\mlimits_{m-n\approx 2^l}u_m u_n\v_{L^2([0,1]\times \mb{R})}\\
	&\le \mlimits_{l\in \mb{N}}\v \mlimits_{m-n\approx 2^l}u_m u_n\v_{L^2([0,1]\times \mb{R})}
\end{align*}

\begin{claim}\label{claim3.10}
For all $ l\in \N$, we have that 
\begineq
\v \mlimits_{m-n\approx 2^l}u_m u_n\v_{L^2([0,1]\times \mb{R})}\lesim (\mlimits_{k\in I, k\in \mb{Z}} \v u_{0,k}\v_{L^2(\mb{R})}^4)^{1/2}.
\endeq
\end{claim}

\noin Notice that $l$ is summed for $\ln|I|$ times, estimate in \eqref{EE3.7} follows directly from this claim.\\

\noin {\bf Proof of Claim \ref{claim3.10}:} Case $l=0$ (i.e. $m\sim n$): 
\begin{align*}
& \v \mlimits_{n\in I, n\in \mb{Z}}u_n^2\v_{L^2([0,1]\times \mb{R})}\lesim (\mlimits_{n\in I, n\in \mb{Z}}\v u_n^2\v_{L^2([0,1]\times \mb{R})}^2)^{1/2}\\
	&\lesim (\mlimits_{n\in I, n\in \mb{Z}} \v u_n\v_{L^4([0,1]\times \mb{R})}^4)^{1/2} \lesim (\mlimits_{n\in I, n\in \mb{Z}} \v u_{0,n}\v_{L^2(\mb{R})}^4)^{1/2},
\end{align*}
where the last step follows from the Strichartz estimate.\\

 Case $l>0$: We split the sum 
\begineq
\sum_{n\in I, n\in \Z}\sum_{m\in I, m\in \Z, m-n\approx 2^l}
\endeq
into
\begineq
\sum_{j\in \Z}\sum_{n\in I, n\approx j 2^l}\sum_{m\in I, m-n\approx 2^l},
\endeq
where $j$ is chosen such that $j2^l, (j+1)2^l\in I$. Hence for $u_n$ with $n\approx j 2^l$ and $u_m$ with $m-n\approx 2^l$, we have that the frequency of the function $u_m u_n$ will be close to $(2j+1) 2^l$, which implies by orthogonality that
\begineq
\begin{split}
& \v \mlimits_{n\in I, n\in \mb{Z}} \mlimits_{m\in I,m\in \mb{Z}, m-n\approx 2^l} u_n u_m\v_{L^2([0,1]\times \mb{R})}\\
&\lesim (\mlimits_{j\in \mb{Z}} \v \mlimits_{n\in I, n\approx j 2^l}\mlimits_{m\in I, m-n\approx 2^l}u_n u_m\v_{L^2([0,1]\times \mb{R})}^2)^{1/2}.
\end{split}
\endeq
By the bilinear estimate, this can be further estimated by
\begin{eqnarray*}
	& & (\mlimits_{j\in \mb{Z}} 2^{-l} \v \mlimits_{n\approx j 2^l}u_{0,n}\v_{L^2(\mb{R})}^2 \v \mlimits_{m\approx (j+1) 2^l}u_{0,m}\v_{L^2(\mb{R})}^2)^{1/2}\\
	&\lesim& (\mlimits_{j\in \mb{Z}} 2^{-l} (\mlimits_{n\approx j 2^l}\v u_{0,n}\v_{L^2(\mb{R})}^2)(\mlimits_{m\approx (j+1) 2^l}\v u_{0,m}\v_{L^2(\mb{R})}^2)^{1/2}\\
	&\lesim& (\mlimits_{j\in \mb{Z}} (\mlimits_{n\approx j 2^l} \v u_{0,n}\v_{L^2(\mb{R})}^4)^{1/2}(\mlimits_{m\approx (j+1) 2^l} \v u_{0,m}\v_{L^2(\mb{R})}^4)^{1/2})^{1/2}\\
	&\lesim& (\mlimits_{n\in I, n\in \mb{Z}}\v u_{0,n}\v_{L^2(\mb{R})}^4)^{1/4}(\mlimits_{m\in I, m\in \mb{Z}}\v u_{0,m}\v_{L^2(\mb{R})}^4)^{1/4}.
\end{eqnarray*}

\noin where the first step is again by the orthogonality in $L^2$,  and the last two steps are just by H\"older's inequality. Hence we have finished the proof of the claim thus the proof of \eqref{EE3.7}.\\

\indent For \eqref{EE3.8}, the starting point is still the same:\\

\noin 
\begineq
\v P_I u \v_{L^4([0,1]\times \mb{R})}^2\lesim \mlimits_{l\in \mb{N}}\v \mlimits_{m-n\approx 2^l}u_m u_n\v_{L^2([0,1]\times \mb{R})}.
\endeq
\noin Case $l=0$ is also the same. For the case $l>0$,\\

\noin $\v \mlimits_{n\in I, n\in \mb{Z}} \mlimits_{m\in I, m-n\approx 2^l} u_n u_m\v_{L^2([0,1]\times \mb{R})}
\lesim (\mlimits_{j\in \mb{N}}\v \mlimits_{n\approx j 2^l}u_n \mlimits_{m\approx (j+1) 2^l} u_m\v_{L^2}^2)^{1/2}$\\

\noin Now we apply the interpolation argument in Proposition \ref{prop2.5}, for the terms $u_{j,l}:=\sum_{n\approx j2^l}u_n$, we can write as a sum $u_{j,l}=u_{1,j,l}+u_{2,j,l}$ with the estimate:\\
\begineq\label{lemma42inter}
\frac{1}{(\ln|I|)^{\beta}}\v u_{1,j,l}\v_{U_{\del}^2}+e^{\epsilon (\ln|I|)^{\beta}}\v u_{2,j,l}\v_{U_{\del}^8}\lesim \v u_{j,l}\v_{V_{\del}^2}.
\endeq

\noin When continuing the estimate, there will be four terms showing up (as we do interpolation for two terms). For the term containing $u_{1, j,l}$ and $u_{1,j+1,l}$, which will be denoted as $\Lambda_1$ (later we still use the notation $u_{j,l}$ defined above for simplicity):
\begin{eqnarray*}
\Lambda_1&\lesim& (\mlimits_{j\in \mb{N}}\v u_{1,j,l}u_{1,j+1,l}\v_{L^2}^2)^{1/2}\\
	&\lesim& (\mlimits_{j\in \mb{N}}2^{-l}\v u_{1,j,l}\v_{U_{\del}^2}^2\v u_{1,j+1,l}\v_{U_{\del}^2}^2)^{1/2}\\
	&\lesim& (\ln|I|)^{2\beta}(\mlimits_{j\in \mb{N}}2^{-l}\v u_{j,l}\v_{V_{\del}^2}^2\v u_{j+1,l}\v_{V_{\del}^2}^2)^{1/2}\\
	&\lesim& (\ln|I|)^{2\beta}(\mlimits_{n\in I, n\in \mb{Z}} \v u_n\v_{V_{\del}^2}^4)^{1/2}
\end{eqnarray*}
where the second step is just by the bilinear estimate, the third is an application of \eqref{lemma42inter}, and the last step is the same as the last summation process in the proof of  \eqref{EE3.7}.\\
\indent Now what is left is to handle the rest three terms in \eqref{lemma42inter}, for which we will do in a uniform way (let us take the term containing only $u_{2,j,l}$ and $u_{2,j+1,l}$ for example, and denote it as $\Lambda_2$):

\begin{eqnarray*}
\Lambda_2 	&=& (\mlimits_{j\in \mb{N}}\v u_{2,j,l} u_{2,j+1,l}\v_{L^2}^2)^{1/2}\\
	&\lesim& (\mlimits_{j\in \mb{N}}\v u_{2,j,l}\v_{L^4}^2 \v u_{2,j+1,l}\v_{L^4}^2)^{1/2}\\
	&\lesim& (\mlimits_{j\in \mb{N}}\v u_{2,j,l}\v_{U_{\del}^8}^2 \v u_{2,j+1,l}\v_{U_{\del}^8}^2)^{1/2}\\
	&\lesim& e^{-2\epsilon (\ln|I|)^{\beta}} 2^{l/2} (\mlimits_{j\in \mb{N}}2^{-l}\v u_{j,l}\v_{V_{\del}^2}^2 \v u_{j+1,l}\v_{V_{\del}^2}^2)^{1/2}\\
	&\lesim& (\mlimits_{n\in I, n\in \mb{Z}} \v u_n\v_{V_{\del}^2}^4)^{1/2}
\end{eqnarray*}
where the second step is by H\"older's inequality, the third step is the Strichartz estimate, the fourth is by applying \eqref{lemma42inter}, the last step just by noticing that $e^{-2\epsilon (\ln|I|)^{\beta}} 2^{l/2}\lesim 1$ due to the fact that $l\le \ln|I|$ and $\beta>1$.\\
\indent Then the last step, summation with respect to $l$ gives rise to the factor $(\ln|I|)^{2\beta+1}$, ending the proof of \eqref{EE3.8}.\\

For \eqref{EE3.9}, what is different from \eqref{EE3.8} is the interpolation argument, here we want to gain a time factor in the tri-linear estimate, which makes the proof somehow technical. The first step is the same as before,
$$\v P_I u \v_{L^4([0,T]\times \mb{R})}^2\lesim \mlimits_{l\in \mb{N}}\v \mlimits_{m-n\approx 2^l}u_m u_n\v_{L^2([0,T]\times \mb{R})},$$

\noin then for the case $l=0$,
\begin{align*}
& \v \mlimits_{n\in I, n\in \mb{Z}}u_n^2\v_{L^2([0,T]\times \mb{R})} \lesim (\mlimits_{n\in I, n\in \mb{Z}}\v u_n^2\v_{L^2([0,T]\times \mb{R})}^2)^{1/2}\\
	&\lesim (\mlimits_{n\in I, n\in \mb{Z}} \v u_n\v_{L^4([0,T]\times \mb{R})}^4)^{1/2} \\
& \lesim T^{1/4}(\mlimits_{n\in I, n\in \mb{Z}} \v u_n\v_{L^8L^4([0,T]\times \mb{R})}^4)^{1/2}\\
	&\lesim T^{1/4}(\mlimits_{n\in I, n\in \mb{Z}} \v u_n\v_{V_{\del}^2}^4)^{1/2},
\end{align*}

\noin where in the last step the Strichartz estimate is applied.\\
\indent For the case $l>0$, still as before apply the orthogonality in $L^2$,
\begineq\label{inter3}
\v \mlimits_{m-n\approx 2^l}u_m u_n\v_{L^2([0,T]\times \mb{R})}\lesim (\mlimits_{j\in \mb{N}}\v \mlimits_{n\approx j 2^l}u_n \mlimits_{m\approx (j+1) 2^l} u_m\v_{L^2}^2)^{1/2}
\endeq

\noin now we can see that the above form is suitable for carrying out the interpolation: for the terms $u_{j,l}:=\sum_{n\approx j2^l}u_n$, we can split it into two terms $u_{j,l}=u_{1,j,l}+u_{2,j,l}$ with the estimate,

\begineq\label{lemmainter2}
\frac{1}{|I|^{\beta}}\v u_{1,j,l}\v_{U_{\del}^2}+e^{\epsilon |I|^{\beta}}\v u_{2,j,l}\v_{U_{\del}^8}\lesim \v u_{j,l}\v_{V_{\del}^2}
\endeq
then use the same notation $\Lambda_1, \Lambda_2$ as before,
\begin{eqnarray*}
\Lambda_1&\lesim& (\mlimits_{j\in \mb{N}}\v u_{1,j,l}u_{1,j+1,l}\v_{L^2}^{2\theta}\v u_{1,j,l}u_{1,j+1,l}\v_{L^2}^{2-2\theta})^{1/2}\\
	&\lesim& (\mlimits_{j\in \mb{N}}\v u_{1,j,l}u_{1,j+1,l}\v_{L^2}^{2\theta}\v u_{1,j,l}\v_{L^4}^{2-2\theta} \v u_{1,j+1,l}\v_{L^4}^{2-2\theta})^{1/2}\\
	&\lesim& T^{\frac{1-\theta}{4}}(\mlimits_{j\in \mb{N}}2^{-l\theta}\v u_{1,j,l}\v_{U_{\del}^2}^2 \v u_{1,j+1,l}\v_{U_{\del}^2}^2)^{1/2}\\
	&\lesim& T^{\frac{1-\theta}{4}}|I|^{2\beta}(\mlimits_{j\in \mb{N}}2^{-l\theta}\v u_{j,l}\v_{V_{\del}^2}^2 \v u_{j+1,l}\v_{V_{\del}^2}^2)^{1/2}\\
	&\lesim& T^{\frac{1-\theta}{4}}|I|^{2\beta+\frac{1-\theta}{2}}(\mlimits_{j\in \mb{N}}2^{-l}\v u_{j,l}\v_{V_{\del}^2}^2 \v u_{j+1,l}\v_{V_{\del}^2}^2)^{1/2}\\
	&\lesim& T^{\frac{1-\theta}{4}}|I|^{2\beta+\frac{1-\theta}{2}}(\mlimits_{n\in I, n\in \mb{Z}} \v u_n\v_{V_{\del}^2}^4)^{1/2}
\end{eqnarray*}

\noin where the second inequality is just by H\"older's inequality, the third inequality is by the bilinear estimate and the Strichartz estimate, the fourth is an application of \eqref{lemmainter2}, and the last summation process is the same as the last summation process in the case $l>0$ in \eqref{EE3.7}.
\begin{eqnarray*}
\Lambda_2 	&=& (\mlimits_{j\in \mb{N}}\v u_{2,j,l} u_{2,j+1,l}\v_{L^2}^2)^{1/2}\\
		&\lesim& (\mlimits_{j\in \mb{N}}\v u_{2,j,l}\v_{L^4}^2 \v u_{2,j+1,l}\v_{L^4}^2)^{1/2}\\
		&\lesim& T^{1/4} (\mlimits_{j\in \mb{N}}\v u_{2,j,l}\v_{U_{\del}^8}^2 \v u_{2,j+1,l}\v_{U_{\del}^8}^2)^{1/2}\\
		&\lesim& T^{1/4} e^{-2\epsilon |I|^{\beta}} 2^{l/2} (\mlimits_{j\in \mb{N}}2^{-l}\v u_{j,l}\v_{V_{\del}^2}^2 \v u_{j+1,l}\v_{V_{\del}^2}^2)^{1/2}\\
		&\lesim& T^{1/4} (\mlimits_{n\in I, n\in \mb{Z}} \v u_n\v_{V_{\del}^2}^4)^{1/2}
\end{eqnarray*}

\noin where the second step is by H\"older's inequality, the third step is the Strichartz estimate and H\"older's inequality in the time variable, the fourth step is an application of \eqref{lemmainter2}, and the last step is the same as the last step in $\Lambda_1$.\\
\indent Then we still need to sum over $l\in \mb{N}$, with $2^l\le |I|$, and this just contributes another factor $\ln|I|$, which can be cancelled if we choose a slightly larger $\beta$. So far we have finished the proof of this lemma. $\square$

\section{Trilinear estimate and proof of Theorem \ref{thm1.4}}

The trilinear estimate below is done on the time interval $[0,T]$, while for simplicity, we still use the notation $U_{\del}^2$, $V_{\del}^2$ instead of $U_{\del,[0,T]}^2$, $V_{\del,[0,T]}^2$.

\begin{thm}\label{thm4.1}
Denote $X_p$ as the function space that we will work in, with the norm given by $\v u\v_{X_p}:=\v \v P_n(u)\v_{U_{\del}^2}\v_{l^p}$, and another space $Y_{q}$ with norm $\v u\v_{Y_q}:= \v \v P_n(u)\v_{V_{\del}^2}\v_{l^q}$. We have the following estimates
\begineq\label{E:4.1}
\v \int_0^t e^{i(t-s)\del}(|u|^2 u)(s)ds \v_{X_p}
\lesim A(T) \v u\v_{X_p}^3,
\endeq
\begineq\label{E:4.2}
\v \int_0^t e^{i(t-s)\del}(u \bar{v} w)(s)ds \v_{X_p}
\lesim A(T) \v u \v_{X_p} \v v \v_{X_p} \v w\v_{X_p}.
\endeq
where the coefficient $A(T)$ depending on the time span $T$ is given by $A(T):=T^{1/2}+T^{1/4}+T^{1/p^{+}}$.
\end{thm}

\begin{rem}
We do not have a coefficient like $A(T)=T^{\alpha}$ essentially because we do not have a good scaling in the function space $M_{2,p}$.
\end{rem}

The argument that Theorem \ref{thm4.1} implies the local well-posedness of the equation \eqref{main} is quite standard, hence we will leave it out.\\

\noin {\bf Proof of Theorem \ref{thm4.1}:} Since the proof of \eqref{E:4.1} and \eqref{E:4.2} are similar, we only write down the details for \eqref{E:4.1}. By the duality relation given in Proposition \ref{prop2.7} and Remark \ref{rem2.10} , it's equivalent to prove the estimate:
$$|\int_{[0,T]\times \mb{R}}u \bar{u}u\bar{v}dxdt|\lesim A(T)\v u\v_{X_p}^3 \v v \v_{Y_{p'}}$$

Now we want to apply the frequency decomposition and write 
\begineq
|\mint_{[0,T]\times \mb{R}}u \bar{u}u\bar{v}dxdt|
\endeq
as 
\begineq
|\mlimits_{\lam_1, \lam_2, \lam_3, n\in \mb{Z}}\int_{[0,T]\times \mb{R}}u_{\lam_1}\bar{u}_{\lam_2}u_{\lam_3}\bar{v}_n dxdt|.
\endeq
In order for the term 
\begineq
\int u_{\lam_1}\bar{u}_{\lam_2}u_{\lam_3}\bar{v}_ndxdt
\endeq
not to vanish, we need the relation for different frequencies
\begineq
\lam_1-\lam_2+\lam_3-n \sim 0.
\endeq

In order to see the interaction of different frequencies, we will divide them into several groups. By symmetry, it's enough to consider the cases where $n$ is the smallest or the second smallest(otherwise consider $-\lam_1, -\lam_2, -\lam_3$ and $-n$), then case 1 would be that the four frequencies are comparable, i.e. 
\begineq
n\sim \lam_1\sim \lam_2\sim \lam_3.
\endeq
The next case, case 2 is when the two smallest frequencies are comparible, i.e. 
\begineq
n\sim \lam_3\ll \lam_1\sim \lam_2.
\endeq 
Then what is left is 
\begineq
n\ll \lam_1\ll \lam_3\ll \lam_2 \text{ or } n\ll \lam_1\sim \lam_3\ll \lam_2,
\endeq 
which we call case 3, or 
\begineq
\lam_1\ll n\ll \lam_2\ll \lam_3 \text{ or } \lam_1\ll n\sim \lam_2\ll \lam_3,
\endeq
which we call case 4.\\

\noin {\bf Case 1:} $\lam_1\sim \lam_2\sim \lam_3\sim n$
\begin{align*}
&\sum_{\lam_1\sim \lam_2\sim \lam_3\sim n}|\int_{[0,T]\times \mb{R}}u_{\lam_1}\bar{u}_{\lam_2}u_{\lam_3}\bar{v}_ndxdt|\\
	&\lesim \mlimits_{\lam_1\sim \lam_2\sim \lam_3\sim n} \v u_{\lam_1}\v_{L^4L^4}\v u_{\lam_2}\v_{L^4L^4}\v u_{\lam_3}\v_{L^4L^4} \v v_n\v_{L^4L^4}\\ 
	&\lesim T^{1/2}\mlimits_{\lam_1\sim \lam_2\sim \lam_3\sim n} \v u_{\lam_1}\v_{L^8L^4}\v u_{\lam_2}\v_{L^8L^4}\v u_{\lam_3}\v_{L^8L^4} \v v_n\v_{L^8L^4}\\
	&\lesim T^{1/2}\mlimits_{\lam_1\sim \lam_2\sim \lam_3\sim n} \v u_{\lam_1}\v_{U_{\del}^2}\v u_{\lam_2}\v_{U_{\del}^2}\v u_{\lam_3}\v_{U_{\del}^2} \v v_n\v_{V_{\del}^2}\\
	&\lesim T^{1/2}\v \v u_{\lam_1}\v_{U_{\del}^2}\v_{l^{3p}}\v \v u_{\lam_2}\v_{U_{\del}^2}\v_{l^{3p}}\v \v u_{\lam_3}\v_{U_{\del}^2}\v_{l^{3p}} \v\v v_n\v_{V_{\del}^2}\v_{l^{p'}}
\end{align*}

\noin where the first two steps are just by H\"older's inequality, the third is the Strichartz estimate, and the fourth is by applying H\"older's inequality with index $1/3p$, $1/3p$, $1/3p$ and $1/p'$.\\

\noin {\bf Case 2:} $n\sim \lam_3\ll \lam_1 \sim \lam_2$.
\begin{align*}
&\mlimits_{n\sim \lam_3\ll \lam_1 \sim \lam_2}|\mint_{[0,T]\times \mb{R}}u_{\lam_1}\bar{u}_{\lam_2} u_{\lam_3}\bar{v}_ndxdt|\\
	&\lesim \mlimits_{n\sim \lam_3\ll \lam_1 \sim \lam_2} \v u_{\lam_1} u_{\lam_3}\v_{L^2L^2}\v u_{\lam_2} v_n\v_{L^2L^2}\\
	&\lesim \mlimits_{\lam_3\ll \lam_1} T^{\theta/2} \frac{\ln^{4(1-\theta)}(\lam_1-\lam_3)}{(\lam_1-\lam_3)^{1-\theta}}\v u_{\lam_1}\v_{V_{\del}^2}^2 \v u_{\lam_3}\v_{V_{\del}^2}\v v_{\lam_3}\v_{V_{\del}^2}\\
	&\lesim T^{\theta/2} \v \mlimits_{\lam_1}\frac{\ln^{4(1-\theta)}(\lam_1-\lam_3)}{(\lam_1-\lam_3)^{1-\theta}}\v u_{\lam_1}\v_{V_{\del}^2}^2 \v_{l^{\infty}_{\lam_3}} \v \v u_{\lam_3}\v_{V_{\del}^2}\v_{l_{\lam_3}^p}\v\v v_{\lam_3}\v_{V_{\del}^2}\v_{l_{\lam_3}^{p'}}\\
	&\lesim T^{\theta/2} \v u\v_{X_p}^3 \v v\v_{Y_{p'}}
\end{align*}

\noin where in the above the first inequality is just by H\"older's inequality, in the second inequality we applied Corollary \ref{cor3.5} and identify $\lam_1$ with $\lam_2$ and $n$ with $\lam_3$ in the summation, the third is again by H\"older's inequality for $\lam_3$ with exponents $l^{\infty}_{\lam_3}$, $l_{\lam_3}^p$ and $l_{\lam_3}^{p'}$, and the last follows form Young's convolution inequality provided that $1-\theta>1-2/p$, namely $\theta/2=1/{p^+}$.\\

\noin {\bf Case 3:} $n\ll \lam_1 \ll \lam_3\ll \lam_2 $ or $n\ll \lam_1 \sim \lam_3\ll \lam_2 $\\
In order to see more orthogonality between different frequencies, we do dyadic decompositions for 
$u_{\lam_1}$,$u_{\lam_2}$ and $u_{\lam_3}$, and keep using uniform decomposition for $v_n$, that is to say, we rewrite the term 
$$\sum_{\lam_1, \lam_2, \lam_3 \in \mb{Z}}\int u_{\lam_1}\bar{u}_{\lam_2}u_{\lam_3}\bar{v}_n dxdt$$
as 
$$\sum_{n\in \mb{Z}, j_1\le j_2 \in \mb{N}}\int \bar{v}_n u_{n+I_{j_1}}u_{n+I_{j_2}}\bar{u}_{n+I_{j_1,j_2}}dxdt,$$
where for $j\in \mb{N}$, by $I_j$ we mean the frequency interval near $2^j$, to be precise,
$$I_j :=[\frac{3}{4}\cdot 2^j, \frac{3}{2}\cdot 2^j].$$
That is to say, we group $\lam_1$ together to the shifted dyadic interval $n+I_{j_1}$, and $\lam_3$ to the shifted dyadic interval $n+I_{j_2}$, then $\lam_2\in n+ I_{j_1, j_2}$ will naturally be grouped by the frequency relation $\lam_1-\lam_2+\lam_3-n\sim 0$, i.e. 
$$I_{j_1,j_2}:=I_{j_1}+I_{j_2}.$$

Let us first calculate the modulation, for $\lam_1\in n+I_{j_1}, \lam_3\in n+I_{j_2}, \lam_2\in n+I_{j_1, j_2}$ satisfying the frequency relation $\lam_1-\lam_2+\lam_3-n\sim 0,$ we get
$$|\lam_1^2-\lam_2^2+\lam_3^2-n^2|\approx |(n-\lam_1)(n-\lam_3)|\approx 2^{j_1} 2^{j_2}.$$

\noin {\bf Case 3.1:} $j_1 \sim j_2$\\

\noin Due to the fact that the high modulation could fall on different terms, we need to discuss this problem in more detail and the techniques would also be different accordingly. Moreover notice that the three larger frequencies are all located near $n+I_{j_1}$ due to the fact that $j_1\sim j_2$, hence they can be somehow identified.\\

\noin {\bf Case 3.1.1} High modulation in either $u_{n+I_{j_1}}$, $u_{n+I_{j_2}}$ or $u_{n+I_{j_1, j_2}}$(in the following calculation we take  $u_{n+I_{j_1}}$ as example):\\

\noin $|\mlimits_{n, j_1, j_2}\int \bar{v}_n u_{n+I_{j_1}} u_{n+I_{j_2}} \bar{u}_{n+I_{j_1, j_2}} dxdt|$\\

\noin $\lesim \mlimits_{n, j_1, j_2} \v v_n\v_{L^{\infty}L^{\infty}} \v u_{n+I_{j_1}}\v_{L^2 L^2} \v u_{n+I_{j_2}}\v_{L^4L^4} \v u_{n+I_{j_1, j_2}}\v_{L^4L^4}$\\


\noin $\lesim \mlimits_{n, j_1, j_2}\v v_n\v_{L^{\infty}L^2} \v u_{n+I_{j_1}}\v_{L^2 L^2} \v u_{n+I_{j_2}}\v_{L^4L^4}^2$\\

\noin $\lesim \mlimits_{n, j_1, j_2} 2^{-\frac{j_1}{2}-\frac{j_2}{2}}(2^{2\beta j_2+\frac{1-\theta}{2} j_2}T^{\frac{1-\theta}{4}}+T^{1/4})...\\
...\v v_n \v_{V_{\del}^2}(\mlimits_{k\in n+I_{j_1}, k\in \mb{Z}}\v u_k\v_{V_{\del}^2}^2)^{1/2}
(\mlimits_{k\in n+I_{j_2}, k\in \mb{Z}}\v u_k\v_{V_{\del}^2}^4)^{1/2}$\\

\noin $\lesim \mlimits_{n,j_1,j_2} 2^{-\frac{j_1}{2}-\frac{j_2}{2}} (2^{2\beta j_2+\frac{1-\theta}{2} j_2}T^{\frac{1-\theta}{4}}+T^{1/4})...\\
... \v v_n \v_{V_{\del}^2} 2^{(1/2-1/p)j_1}(\mlimits_{k\in n+I_{j_1}} \v u_k\v_{U_{\del}^2}^p)^{1/p}
 2^{(1/2-2/p)j_2}(\mlimits_{k\in n+I_{j_2}}\v u_k\v_{U_{\del}^2}^p)^{2/p}$\\

\noin $\lesim \mlimits_{n, j_1, j_2} 
2^{-\frac{j_1+2j_2}{p}+2\beta j_2+\frac{1-\theta}{2} j_2}( T^{\frac{1-\theta}{4}}+T^{1/4})
\v v_n\v_{V_{\del}^2}(\mlimits_{k\in n+I_{j_1}} \v u_k\v_{U_{\del}^2}^p)^{1/p} 
\v u\v_{X_p}^2$\\

\noin $\lesim \mlimits_{n, j_1} 
(T^{\frac{1-\theta}{4}}+T^{1/4}) \v v_n\v_{V_{\del}^2} 2^{-\frac{j_1}{p}-\epsilon j_1} 
(\mlimits_{k\in n+I_{j_1}, k\in \mb{Z}} 
\v u_k\v_{U_{\del}^2}^p)^{1/p}
\v u \v_{X_p}^2$\\

\noin $\lesim (T^{\frac{1-\theta}{4}}+T^{1/4})(\mlimits_n \v v_n\v_{V_{\del}^2}^{p'})^{1/{p'}} 
\mlimits_{j_1} 2^{-\epsilon j_1} (\mlimits_n (2^{-j_1} \mlimits_{k\in n+I_{j_1}, k\in \mb{Z}}\v u_k \v_{U_{\del}^2}^p))^{1/p} \v u \v_{X_p}^2$\\

\noin where the first inequality is by H\"older's inequality, the second is by Bernstein's inequality and identifying $n+I_{j_2}$ with $n+I_{j_1, j_2}$, the third is by applying Lemma \ref{lem3.6}, the fourth is by applying H\"older's inequality to the two latter brackets, the fifth is just a rewriting, the sixth holds if we choose $\theta$, $\beta$ and small enough $\epsilon$ such that $\frac{1-\theta}{2}+2\beta+\epsilon <\frac{2}{p}$, which implies 
$$2^{-\frac{2j_2}{p}+2\beta j_2+\frac{1-\theta}{2}j_2}\lesim 2^{-\epsilon j_1},$$
and the last is by applying H\"{o}lder's inequality with exponents $l_{j_1}^{\infty}l_n^{p'}$ and $l_{j_1}^1l_n^p$, and in $$\mlimits_n 
(2^{-j_1} \mlimits_{k\in n+I_{j_1}, k\in \mb{Z}}\v u_k \v_{U_{\del}^2}^p),$$
every $\v u_k \v_{U_{\del}^2}^p$ is added for $2^{j_1}$ times, which cancels the coefficient in front, giving exactly what we want on the right hand side of Theorem \ref{thm4.1}.\\
\indent For the choice of $\theta$, we observe that as $\beta$ in \eqref{EE3.9} of Lemma \ref{lem3.6} and the small parameter $\epsilon$ above can be chosen arbitrarily small, it's enough to choose 
$$\frac{1-\theta}{2}<\frac{2}{p},$$ 
that is, $\frac{1-\theta}{4}=\frac{1}{p^+}$. Hence we have finished the proof of The Case 3.1.1.\\

\noin {\bf Case 3.1.2} High modulation in $v_n$:
\begin{eqnarray*}
	& & \mlimits_{n, j_1, j_2} |\int \bar{v}_n u_{n+I_{j_1}} u_{n+I_{j_2}} \bar{u}_{n+I_{j_1, j_2}}dxdt|\\
	&\lesim& \mlimits_{n, j_1, j_2} 
\v v_n \v_{L^2L^{\infty}} \v u_{n+I_{j_1}}\v_{L^{\infty}L^2} 
\v u_{n+I_{j_2}}\v_{L^4L^4} \v u_{n+I_{j_1, j_2}}\v_{L^4L^4}\\
	&\lesim& \mlimits_{n, j_1, j_2} 
\v v_n\v_{L^2L^2} \v u_{n+I_{j_1}}\v_{V_{\del}^2} \v u_{n+I_{j_2}}\v_{L^4L^4} \v u_{n+I_{j_1, j_2}}\v_{L^4L^4}
\end{eqnarray*}
where in the second inequality we applied the Berstein estimate for the first term and the Strichartz estimate for the second term. Now the strategy is the same as in the previous case: to gain from high modulation in the term $\v v_n\v_{L^2L^2}$, to apply the $V^2$ orthogonality argument for the term $\v u_{n+I_{j_1}}\v_{V_{\del}^2}$, and to carry exactly the same $L^4L^4$ estimate for the last two terms as before.\\

\noin {\bf Case 3.2:} $j_1 \ll j_2.$\\

\noin {\bf Case 3.2.1} High modulation in $u_{n+I_{j_1}}$:
\begin{eqnarray*}
	& & \mlimits_{n, j_1, j_2} |\int \bar{v}_n u_{n+I_{j_1}} u_{n+I_{j_2}} \bar{u}_{n+I_{j_1, j_2}}dxdt|\\
	&\lesim& \mlimits_{n, j_1, j_2} 
\v v_n \v_{L^{\infty}L^{\infty}}\v u_{n+I_{j_1}}\v_{L^2L^2} 
\v u_{n+I_{j_2}}\v_{L^4L^4} \v u_{n+I_{j_1, j_2}}\v_{L^4L^4}
\end{eqnarray*}

\noin which is then the same as case 3.1.1.\\

\noin {\bf Case 3.2.2} High modulation in $v_n$:
\begin{eqnarray*}
	& & \mlimits_{n, j_1, j_2} |\int \bar{v}_n u_{n+I_{j_1}} u_{n+I_{j_2}} \bar{u}_{n+I_{j_1, j_2}}dxdt|\\
	&\lesim& \mlimits_{n, j_1, j_2} 
\v v_n \v_{L^2L^{\infty}}\v u_{n+I_{j_1}}\v_{L^{\infty}L^2} 
\v u_{n+I_{j_2}}\v_{L^4L^4} \v u_{n+I_{j_1, j_2}}\v_{L^4L^4}
\end{eqnarray*}

\noin which is the same as case 3.1.2.\\

\noin {\bf Case 3.2.3} High modulation in either $u_{n+I_{j_2}}$ or $u_{n+I_{j_1, j_2}}$:\\

\noin As we are in this case $j_1\ll j_2$, it is not difficult to see that the frequency of $u_{n+I_{j_1, j_2}}$ is also localised near $n+I_{j_2}$, which explains why this two term can be identified. Hence w.l.o.g we assume that the high modulation lies in $u_{n+I_{j_1, j_2}}$. In order to gain from high modulation, we need to apply $L^2L^2$ estimate to the term $u_{n+I_{j_1, j_2}}$. For the other terms, we just need to observe the frequency separation of $n+I_{j_1}$ and $n+I_{j_2}$ is $2^{j_2}$ because of $j_1\ll j_2$, which makes it possible to apply one bilinear estimate. To be precise:
\begin{eqnarray*}
	& & \mlimits_{n, j_1, j_2} |\int \bar{v}_n u_{n+I_{j_1}} u_{n+I_{j_2}} \bar{u}_{n+I_{j_1, j_2}}dxdt|\\
	&\lesim& \mlimits_{n, j_1, j_2} 
\v v_n \v_{L^{\infty}L^{\infty}} \v u_{n+I_{j_1}} u_{n+I_{j_2}}\v_{L^2L^2} \v u_{n+I_{j_1, j_2}}\v_{L^2L^2}
\end{eqnarray*}

\noin then the following is basically the same as The Case 3.1.1. Especially for the bilinear estimate: in order to gain a small time factor, we will not apply the full bilinear estimate to the term $\v u_{n+I_{j_1}} u_{n+I_{j_2}}\v_{L^2L^2}$, but the bilinear estimate together with H\"older's inequality as in Case 2 and Case 3.1.1. \\

\noin {\bf Case 4:} $\lam_1\ll n\ll \lam_2\ll \lam_3$ or $\lam_1\ll n\sim \lam_2\ll \lam_3$\\

\noin As explained at the beginning of Case 3, we will also do dyadic decomposition for the frequency $\lam_1$, $\lam_2$ and $\lam_3$, and still keep using uniform decomposition for $n$ , i.e. we will write the term 
\begineq
\sum_{\lam_1, \lam_2, \lam_3 \in \mb{Z}}\int u_{\lam_1}\bar{u}_{\lam_2}u_{\lam_3}\bar{v}_n dxdt
\endeq
as 
\begineq
\sum_{n\in \mb{Z}, j_1, j_2 \in \mb{N}}\int \bar{v}_n u_{n-I_{j_1}}u_{n+I_{j_2}}\bar{u}_{n+I_{-j_1, j_2}}dxdt,
\endeq
where $\lam_1$ is grouped to the shifted dyadic interval $n-I_{j_1}$, $\lam_2$ to the shifted dyadic interval $n+I_{j_2}$, and $\lam_3\in n+I_{-j_1, j_2}$ will again be grouped by the frequency relation $\lam_1-\lam_2+\lam_3-n\sim 0$.

The relation of $j_1$ and $j_2$ will determine how different frequencies are separated, hence we will still divide case 4 into three subcases, namely $j_1\ll j_2$, $j_1\sim j_2$ and $j_1\gg j_2$. However, compared with the subcases in Case 3, the calculation here is still a combination of H\"older's inequality, Bernstein's inequality, the Strichartz estimate, the bilinear estimate, the $L^4$ estimate in Lemma \ref{lem3.6} and the $L^2L^2$ estimate for the term with high modulation. Hence we will leave the details out.\\

\noin So far the proof for all the cases in the trilinear estimate is completed.$\square$

\section{Orlicz Spaces}

\begin{defi}(\cite{Kuf})\label{defi5.1}
We call a $C^2$ convex function $\Phi: \mb{R}^+\to \mb{R}^+$ a Young function if 
$\Phi$ satisfies the following conditions: $\Phi(0)=0; \Phi(s)>0$ for $s>0$; $\lim_{s\to 0^+}\Phi'(s)=0$; $\lim_{s\to \infty}\Phi'(s)=\infty$. 
\end{defi}

\begin{rem}
The Young functions we use in the context are always smooth, which explains why we state the definition only for $C^2$ functions. There is a more general definition which can be found in \cite{Kuf}.
\end{rem}

\begin{defi}\label{defi5.2}(\cite{Kuf})
Let $\Phi$ be a Young function, the function space $L^{\Phi}$ contains all measurable functions on $\R^n$ such that the following norm is finite:
$$\v u\v_{L^{\Phi}}:=\inf\{k>0: \int_{\R^n}\Phi(\frac{|u(x)|}{k})dx\le 1\}.$$
The function space $l^{\Phi}$ contains all the complex-valued sequence $\{u_n\}_{n\in \N}$ such that the discrete version of the above norm is finite:
\begineq
\|\{u_n\}_{n\in \N}\|_{l^{\Phi}}:=\inf\{k>0: \sum_{n\in \N}\Phi(\frac{|u_n|}{k})\le 1\}.
\endeq
\end{defi}

\begin{defi}(\cite{Kuf})
For a Young function $\Phi$, its convex conjugate function $\Psi$ is defined by
$$\Psi(t):=\sup_{s\in\mb{R}^+}\{st-\Phi(s)\}.$$
\end{defi}

\begin{rem}
It is not difficult to see that $\Psi$ is also a Young function.
\end{rem}

\begin{prop}\label{prop5.6}(General H\"older's inequality; \cite{Kuf})
Let $\Phi, \Psi$ be a pair of convex conjugate functions, then
\begineq
\int_{\R^n}|u(x)v(x)|dx \lesim \v u\v_{L^{\Phi}}\v v\v_{L^{\Psi}}.
\endeq
\end{prop}
\noindent {\bf Proof of Proposition \ref{prop5.6}:} Take $u\in L^{\Phi}$ with $\v u\v_{L^{\Phi}}=1$, and $u\in L^{\Psi}$ with $\v u\v_{L^{\Psi}}=1$. Then by the definition of convex conjugate functions, we have
$$\int |uv| dx\le \int \Phi(u)+\Psi(v) dx \le 2,$$
where the second step is by the definition of Orlicz spaces and the normalization we did for $u$ and $v$. $\square$

\begin{prop}(\cite{Kuf})
Under the above notations, for a fixed function $v$ in $L^{\Psi}(\R^n)$, the expression
$$F(u):=\int_{\R^n}u(x)v(x)dx, u\in L^{\Phi}(\R^n)$$
defines a continuous linear functional $F$ on $L^{\Phi}(\R^n)$ and 
$$\frac{1}{2}\v v\v_{L^{\Psi}}\le \v F\v \le \v v\v_{L^{\Psi}}.$$
\end{prop}

In the following we will collect some technical lemmas which will be used in the proof of the trilinear estimate in the Orlicz spaces.

\begin{lem}\label{lem5.8}
With $\bar{\Phi}=e^{-(\frac{1}{x})^{\frac{1}{2 \gamma}}+C_{\gamma}x^{\frac{1}{2}}}$,$\gamma>3/2$, and $\bar{\Psi}$ the corresponding convex conjugate function, we have that the sequence $\{\frac{\ln^2 n}{n}\}_{n\in\mb{N}_+}\in l^{\bar{\Psi}}$, i.e. 
\begineq
\v \{\frac{\ln^2 n}{n}\}_{n\in\mb{N}_+}\v_{l^{\bar{\Psi}}}<\infty.
\endeq
\end{lem}

\noin {\bf Proof of Lemma \ref{lem5.8}:} For the function $g(s):=st-\bar{\Phi}(s)$ in the definition of the convex conjugate function, we first observe that\\
$$\bar{\Phi}(s)\approx e^{-(\frac{1}{s})^{\frac{1}{2\gamma}}}$$
when $s$ small, which further implies that
\begineq\label{E:5.2}
g'(s)\approx t-e^{-(\frac{1}{s})^{\frac{1}{2\gamma}}} \frac{1}{2\gamma}(\frac{1}{s})^{\frac{1}{2\gamma}+1}.
\endeq
Now we solve $g'(s)=0$ by making the Ansatz\\
$$s=C(t)[\ln(\frac{1}{t})]^{-2\gamma},$$
substituting which into \eqref{E:5.2} we obtain
\begineq\label{E:5.3}
2\gamma t\approx t^{(\frac{1}{C(t)})^{\frac{1}{2\gamma}}}(\ln\frac{1}{t})^{2\gamma+1}(\frac{1}{C(t)})^{\frac{2\gamma+1}{2\gamma}},
\endeq
from this expression we can see that $C(t)$ must tend to 1 when $t$ tends to 0.\\
\indent Substitute $s$ into the definition of convex conjugate function and apply \eqref{E:5.3}:
\begin{eqnarray*}
st-\bar{\Phi}(s) &\approx& C(t) t [\ln(\frac{1}{t})]^{-2\gamma}-t^{(\frac{1}{C(t)})^{\frac{1}{2\gamma}}}\\
	&\approx& C(t) t [\ln(\frac{1}{t})]^{-2\gamma}- 2\gamma t (\ln\frac{1}{t})^{-2\gamma-1}(C(t))^{\frac{2\gamma+1}{2\gamma}}\\
	&\approx& C(t) t [\ln(\frac{1}{t})]^{-2\gamma} \approx t [\ln(\frac{1}{t})]^{-2\gamma}
\end{eqnarray*}
when $t$ small.\\
\indent Concerning the norm $\v\{\frac{\ln^2 n}{n}\}\v_{l^{\bar{\Psi}}},$ by the definition of the Orlicz space, it would suffice if we can prove that the sequence $\{\frac{\ln^2 n}{n} (\ln \frac{n}{\ln^2 n})^{-2\gamma}\}$ is summable, while this is the case as
$$\mlimits_n \frac{\ln^2 n}{n} (\ln \frac{n}{\ln^2 n})^{-2\gamma}\approx \mlimits_n \frac{\ln^2 n}{n} (\ln n)^{-2\gamma},$$
where $\gamma>3/2$ will be enough. This finishes the proof of the above lemma.$\square$

\begin{lem}\label{lem5.9}
For $N\in \mb{Z}_+$, for the sequence $\{a_N(n)\}_{n\in \mb{Z}}$ with
\begin{displaymath}
a_N(n)=\left\{ \begin{array}{ll}								
1& {1\le n \le N}\\
0 & {\text{else}}
\end{array} \right.
\end{displaymath}
we have that $$\v \{a_N(n)\}\v_{l^{\Psi_3}}\lesim \frac{N}{(\ln N)^{4\gamma}},$$
with $\Phi_3(x)=e^{-(\frac{1}{x})^{\frac{1}{4\gamma}}+C_{\gamma} x^{\frac{1}{4}}}$ and $\Psi_3$ the convex conjugate function of $\Phi_3$.
\end{lem}

\noin {\bf Proof of Lemma \ref{lem5.9}:} By the proof of Lemma \ref{lem5.8} we have for $t$ small that 
\begineq\label{E:5.4}
\Psi_3(t)\approx t(\ln \frac{1}{t})^{-4\gamma}.
\endeq
By the definition of the Orlicz space, we obtain
$$\v \{a_N(n)\}\v_{l^{\Psi_3}}= \inf \{k>0: \Psi_3(\frac{1}{k})\le \frac{1}{N}\},$$
Applying \eqref{E:5.4}, we get
$$\frac{1}{k}(\ln k)^{4\gamma}\lesim \frac{1}{N},$$
which further implies 
$$k\gtrsim \frac{N}{(\ln N)^{4\gamma}},$$ 
which is a upper bound for the norm. $\square$

\begin{lem}\label{lem5.10}
For any $j\in \mb{Z}_+$, any sequence $\{a_k\}_{k\in \mb{Z}}$ and Young function $\Phi$, the following inequality holds true:
\begineq\label{E:5.5}
\v \{ (\mlimits_{k\in n+I_{2^j}, k\in \mb{Z}} 
2^{-j} a_k^2)^{1/2}\}_{n\in \mb{Z}}\v_{l^{\Phi}}
\lesim \v \{a_k\}_{k\in \mb{Z}} \v_{l^{\Phi}}.
\endeq
\end{lem}

\noin {\bf Proof of Lemma \ref{lem5.10}:} The proof is mainly by applying Jensen's inequality:
\begin{eqnarray*}
	&&\v \{ (\mlimits_{k\in n+I_{2^j}, k\in \mb{Z}} 2^{-j} a_k^2)^{1/2}\}_{n\in \mb{Z}}\v_{l^{\Phi}}\\
	&=& \inf \{\lam>0: \mlimits_n \Phi(\frac{(\sum_{k\in n+I_{2^j}, k\in \mb{Z}} 2^{-j} a_k^2)^{1/2}}{\lam})\le 1\}\\
	&=&\inf \{\lam^{1/2}>0: \mlimits_n \bar{\Phi}(\frac{\sum_{k\in n+I_{2^j}, k\in \mb{Z}} 2^{-j} a_k^2}{\lam})\le 1\}\\
	&\le& \inf \{\lam^{1/2}>0: \mlimits_n \mlimits_{k\in n+I_{2^j}, k\in \mb{Z}} 2^{-j} \bar{\Phi}(\frac{a_k^2}{\lam})\le 1\}\\
	&\le& \inf\{\lam>0: \mlimits_n \mlimits_{k\in n+I_{2^j}, k\in \mb{Z}} 2^{-j} \Phi(\frac{a_k}{\lam})\le 1\}\\
	&\le& \inf\{\lam>0: \mlimits_{k\in \mb{Z}} \Phi(\frac{a_k}{\lam})\le 1\},
\end{eqnarray*}

\noin where the first step is by the definition of the Orlicz space, the second is a rewriting with $\bar{\Phi}$ given by $\bar{\Phi}(t):= \Phi(t^{1/2})$, which is again a Young function, the third is by Jensen's inequality, the fourth is the same as what we did for the second step, and the last step is just by counting how many times every term $\Phi(\frac{a_k}{\lam})$ is added.$\square$

\section{Trilinear Estimate in Orlicz Spaces and Proof of Theorem \ref{thm1.8}}
In this section, we will prove a trilinear estimate in the Orlicz space similar to the one in Section 4. The estimate done below is on the time interval $[0,1]$, for simplicity we use the notation $U_{\del}^2$ and $V_{\del}^2$ instead of $U_{\del,[0,1]}^2$ and $V_{\del,[0,1]}^2$

\begin{thm}\label{thm6.1}(small data well-posedness on the time interval $[0,1]$)
 We have the following estimates
$$\v \int_0^t e^{i(t-s)\del}(|u|^2 u)(s)ds \v_{X_{\Phi}}
\lesim \v u\v_{X_{\Phi}}^3,$$
$$\v \int_0^t e^{i(t-s)\del}(u \bar{v} w)(s)ds \v_{X_{\Phi}}
\lesim \v u\v_{X_{\Phi}} \v v\v_{X_{\Phi}} \v w\v_{X_{\Phi}},$$
where $\Phi$ is the Young function defined in Theorem \ref{thm1.4}, $\Psi$ is the convex conjugate function of $\Phi$, $X_{\Phi}$ is the function space that we will work in, with its norm given by $\v u\v_{X_{\Phi}}:=\v\v P_n(u)\v_{U_{\del}^2}\v_{l_n^{\Phi}}$ and $Y_{\Psi}$ is another space that we need, with a norm $\v u\v_{Y_{\Psi}}:= \v\v P_n(u)\v_{V_{\del}^2}\v_{l_n^{\Psi}}$.
\end{thm}

\noin {\bf Proof of Theorem \ref{thm6.1}}: Again we will only write down the proof of the first trilinear estimate, as the second will follow in a similar manner. By duality, it is enough to prove the estimate:
$$|\int_{[0,1]\times \mb{R}}u \bar{u}u\bar{v}dxdt|\lesim \v u\v_{X_{\Phi}}^3 \v v \v_{Y_{\Psi}}.$$

\noin {\bf Case 1:} $\lam_1\sim \lam_2\sim \lam_3\sim n$
\begin{align*}
& \mlimits_{\lam_1\sim \lam_2\sim \lam_3\sim n}|\int_{[0,1]\times \mb{R}}u_{\lam_1}\bar{u}_{\lam_2}u_{\lam_3}\bar{v}_n dxdt|\\
	&\lesim \mlimits_n \v u_n\v_{L^4L^4}^3 \v v_n\v_{L^4L^4} \lesim \mlimits_n \v u_n\v_{L^8L^4}^3 \v v_n\v_{L^8L^4}\\
	&\lesim \mlimits_n \v u_n\v_{U_{\del}^2}^3 \v v_n\v_{V_{\del}^2} \lesim \v \v u_n\v_{U_{\del}^2}\v_{l^{\infty}}^2\v \v u_n\v_{U_{\del}^2}\v_{l^{\Phi}} \v\v v_n\v_{V_{\del}^2}\v_{l^{\Psi}}\\
	&\lesim \v u_n\v_{X_{\Phi}}^3 \v v_n\v_{Y_{\Psi}}.
\end{align*}

\noin Here the first two inequalities are by H\"older's inequality, the third is by the Strichartz estimate, the fourth is H\"older's inequality for the general Orlicz space, and the last step is by the trivial embedding $l^{\Phi}\emb l^{\infty}$.\\

\noin {\bf Case 2:} $\lam_1 \sim \lam_2 \ll \lam_3 \sim n$.\\
\begin{align*}
&\mlimits_{\lam_1 \sim \lam_2 \ll \lam_3 \sim n}|\mint_{[0,1]\times \mb{R}}u_{\lam_1}\bar{u}_{\lam_2} u_{\lam_3}\bar{v}_n dxdt|\\
	&\le \mlimits_{\lam_1\ll \lam_3} \frac{\ln^2(\lam_3-\lam_1)}{\lam_3-\lam_1} \v u_{\lam_1}\v_{U_{\del}^2}^2 \v u_{\lam_3}\v_{U_{\del}^2}\v v_{\lam_3}\v_{V_{\del}^2}\\
	&\lesim \v \mlimits_{\lam_1} \frac{\ln^2 (\lam_3-\lam_1)}{\lam_3-\lam_1} \v u_{\lam_1}\v_{U_{\del}^2}^2\v_{l^{\infty}_{\lam_3}} 
\v\v u_{\lam_3}\v_{U_{\del}^2}\v_{l^{\Phi}} \v\v v_{\lam_3}\v_{V_{\del}^2}\v_{l^{\Psi}}\\
	&\lesim \v \{\frac{\ln^2 n}{n}\}\v_{l^{\bar{\Psi}}}\v \v u_{\lam_1}\v_{U_{\del}^2}^2\v_{l^{\bar{\Phi}}}
\v\v u_{\lam_3}\v_{U_{\del}^2}\v_{l^{\Phi}} \v\v v_{\lam_3}\v_{V_{\del}^2}\v_{l^{\Psi}}\\
	&\lesim \v \{\frac{\ln^2 n}{n}\}\v_{l^{\bar{\Psi}}}\v \v u_{\lam_1}\v_{U_{\del}^2}\v_{l^{\Phi}}^2
\v\v u_{\lam_3}\v_{U_{\del}^2}\v_{l^{\Phi}} \v\v v_{\lam_3}\v_{V_{\del}^2}\v_{l^{\Psi}}\\
	&\lesim \v u\v_{X_{\Phi}}^3 \v v\v_{Y_{\Psi}}
\end{align*}

\noin where $\bar{\Phi}=e^{-(\frac{1}{x})^{\frac{1}{2\gamma}}+C_{\gamma}x^{\frac{1}{2}}}$, and $\bar{\Psi}$ is the corresponding convex conjugate function. The first inequality is an application of Corollary \ref{cor3.4} and Corollary \ref{cor3.5}, the second and the third are simply by H\"older's inequality for the general Orlicz spaces, the fourth is by the definition of the Orlicz spaces, and the last is by Lemma \ref{lem5.8}.\\

\noin {\bf Case 3:} $n\ll \lam_1 \ll \lam_3 \ll \lam_2$ and $n\ll \lam_1 \sim \lam_3 \ll \lam_2$\\

For all the subcases of this case, the basic strategy is the same as the one used in the proof of Theorem \ref{thm4.1}, the only difference is that in the summation process we need to apply H\"{o}lder's inequality for the general Orlicz space. Hence in the following we will just detail case 3.1.1 for simplicity. Moreover we will leave out the proof for case 4.\\

\noin {\bf Case 3.1.1} high modulation in $u_{n+I_{j_1}}$: \\

\noin $\mlimits_{n, j_1, j_2}|\mint \bar{v}_n u_{n+I_{j_1}} u_{n+I_{j_2}} \bar{u}_{n+I_{j_1, j_2}} dxdt|$\\

\noin $\lesim \mlimits_{n, j_1, j_2} \v v_n\v_{L^{\infty}L^{\infty}} \v u_{n+I_{j_1}}\v_{L^2L^2} \v u_{n+I_{j_2}}\v_{L^4L^4}\v u_{n+I_{j_1, j_2}}\v_{L^4L^4}$\\

\noin $\lesim \mlimits_{n, j_1, j_2}\v v_n\v_{L^{\infty}L^2} \v u_{n+I_{j_1}}\v_{L^2 L^2} \v u_{n+I_{j_2}}\v_{L^4L^4}^2$\\

\noin $\lesim \mlimits_{n, j_1, j_2} 2^{-j_1/2-j_2/2}j_2^{2 \beta+1}\v v_n\v_{V_{\del}^2}(\mlimits_{k\in n+I_{j_1}, k\in \mb{Z}}\v u_{k}\v_{V_{\del}^2}^2)^{1/2} (\mlimits_{k\in n+I_{j_2}, k\in\mb{Z}} \v u_k\v_{V_{\del}^2}^4)^{1/2}$\\

\noin $\lesim \mlimits_{n, j_1, j_2} 
2^{-\frac{j_1}{2}-\frac{j_2}{2}} j_2^{2\beta+1} \v v_n \v_{V_{\del}^2} 
(\mlimits_{k\in n+I_{j_1}, k\in \mb{Z}} 
\v u_k\v_{U_{\del}^2}^2)^{1/2} 
\v \{a_{2^{j_2}}(n)\}\v_{l^{\Psi_3}}^{1/2} \v \v u_k\v_{U_{\del}^2}^4 \v_{l^{\Phi_3}}^{1/2}$

\noin $\lesim \mlimits_{n, j_1} 2^{-\frac{j_1}{2}} j_1^{-(2\gamma-2\beta-1)} \v v_n \v_{V_{\del}^2} 
(\mlimits_{k\in n+2^{j}, k\in \mb{Z}} 
\v u_k\v_{U_{\del}^2}^2)^{1/2} 
\v \v u_k\v_{U_{\del}^2} \v_{l^{\Phi}}^{2}$\\

\noin $\lesim \v \v v_n \v_{V_{\del}^2} \v_{l^{\Psi}} 
(\mlimits_{j_1\in \bf{N^+}} j_1^{-(2\gamma-2\beta-1)}) \v \{ (\mlimits_{k\in n+2^{j_1}, k\in \mb{Z}} 
2^{-j_1} \v u_k\v_{U_{\del}^2}^2)^{1/2}\}_{n\in \mb{Z}}\v_{l^{\Phi}}
\v \v u_k\v_{U_{\del}^2} \v_{l^{\Phi}}^2$\\

\noin where the first step is by H\"older's inequality, the second is by Bernstein's inequality, the third is by applying \eqref{EE3.8} in Lemma \ref{lem3.6}, the fourth is by applying H\"{o}lder's inequality Proposition \ref{prop5.6} to the two brackets, the second last is by Lemma \ref{lem5.9} with the same $\Phi_3$ and $\Psi_3$, and the last is by applying H\"{o}lder's inequality with index $l_{j_1}^{\infty}l_n ^{\Psi}$ and $l_{j_1}^1 l_n ^{\Phi}$. In the end, Lemma \ref{lem5.10} concludes the proof for this case.\\
\indent Now let us see how to choose $\gamma$. In order for the term $\sum_{j_1\in \mb{N^+}} j_1^{-(2\gamma-2\beta-1)}$ in the last inequality to be finite, we need $\gamma>\beta+1$. As $\beta$ can be chosen arbitrarily close to 1 in \eqref{EE3.8} of Lemma \ref{lem3.6}, we obtain the restriction that $\gamma>2$.\\ 
\indent As stated above, we will leave out the proof for the rest cases, hence proof of Theorem \ref{thm6.1} is finished.$\square$\\

Now we turn to the proof of Corollary \ref{cor1.13}, which gives the almost global well-posedness.\\

\noin {\bf Proof of Corollary \ref{cor1.13}:} Following from Theorem \ref{thm1.8} and the embedding \eqref{E:1.5}, well-posedness on the time interval $[0,1]$ is guaranteed. Then it will be enough if we can prove the following scaling: 
\begineq\label{E:6.1}
\v v_0\v_{\hat{L}^{\Phi}}\lesim (\v u_0\v_{\hat{L}^{\Phi}}+\v \hat{u}_0\v_{L^{\infty}})(\ln \lam)^{\gamma}
\endeq
with $v_0(x):=\lam u_0(\lam x)$, and $\lam$ large. This is because for initial data $u_0$ with $\v u_0\v_{\hat{L}^{\Phi}}+\v \hat{u}_0\v_{L^{\infty}}$ being small, we can choose $\lam$ such that 
$$\v v_0\v_{l^{\Phi}L^2}\lesim \v v_0\v_{\hat{L}^{\Phi}}\approx \epsilon_0,$$
where $\epsilon_0$ is the small parameter for Theorem \ref{thm1.8} to hold, then there exists solution $v(t,x)$ with initial data $v_0$ on time interval $[0,1]$. By scaling, 
$$u(t,x):=\frac{1}{\lam}v(\frac{t}{\lam^2},\frac{x}{\lam})$$ 
also solves the same equation with initial data $u_0$ for $0\le t\le \lam^2$.\\

\indent {\bf Proof of \eqref{E:6.1}}: For simplicity, we denote $\v \hat{u}_0\v_{L^{\infty}}=a$,  $\v u_0\v_{\hat{L}^{\Phi}}=b$, and assume w.l.o.g. that $a+b=1$. Then by the definition of $\hat{L}^{\Phi}$, we have that:\\
\begineq\label{defi}
\int_{\mb{R}} e^{-(\frac{b}{|\hat{u}_0(\xi)|})^{1/{\gamma}}+C_{\gamma}\frac{|\hat{u}_0(\xi)|}{b}}d\xi=1.
\endeq
Now we calculate the norm of $v_0$:
\begin{eqnarray*}
\v v_0\v_{\hat{L}^{\Phi}}	&=&	\inf\{M>0| \int_{\mb{R}} e^{-(\frac{M}{|\hat{v}_0(\xi)|})^{1/{\gamma}}+C_{\gamma}\frac{|\hat{v}_0(\xi)|}{M}}d\xi\le 1\}\\
				&=&	\inf\{M>0| \int_{\mb{R}} e^{-(\frac{M}{|\hat{u}_0(\xi)|})^{1/{\gamma}}+C_{\gamma}\frac{|\hat{u}_0(\xi)|}{M}}d\xi\le 1/{\lam}\}
\end{eqnarray*}
Hence it suffices to prove that  there exists $C>0$ such that:
$$\int_{\mb{R}} e^{-(\frac{C (\ln \lam)^{\gamma}}{|\hat{u}_0(\xi)|})^{1/{\gamma}}+C_{\gamma}\frac{|\hat{u}_0(\xi)|}{C (\ln \lam)^{\gamma}}}d\xi\le 1/{\lam}, \any \lam \ge e,$$
or equivalently (replacing $\ln \lam$ by $\lam$):
$$\int_{\mb{R}} e^{\lam-(\frac{C \lam^{\gamma}}{|\hat{u}_0(\xi)|})^{1/{\gamma}}+C_{\gamma}\frac{|\hat{u}_0(\xi)|}{C \lam^{\gamma}}}d\xi\le 1, \any \lam \ge 1.$$
Let us first check for $\lam=1$:

\begineq\label{inte}
\int_{\mb{R}} e^{1-(\frac{C}{|\hat{u}_0(\xi)|})^{1/{\gamma}}+C_{\gamma}\frac{|\hat{u}_0(\xi)|}{C}}d\xi
\endeq

\begineq
=\int_{\mb{R}}e^{-(\frac{b}{|\hat{u}_0(\xi)|})^{1/{\gamma}}+C_{\gamma}\frac{|\hat{u}_0(\xi)|}{b}+(\frac{b}{|\hat{u}_0(\xi)|})^{1/{\gamma}}-C_{\gamma}\frac{|\hat{u}_0(\xi)|}{b}+1-(\frac{C}{|\hat{u}_0(\xi)|})^{1/{\gamma}}+C_{\gamma}\frac{|\hat{u}_0(\xi)|}{C}}
\endeq

\noin Noticing that \eqref{defi} holds true, hence it will be enough to prove that
$$(\frac{b}{|\hat{u}_0(\xi)|})^{1/{\gamma}}-C_{\gamma}\frac{|\hat{u}_0(\xi)|}{b}+1-(\frac{C}{|\hat{u}_0(\xi)|})^{1/{\gamma}}+C_{\gamma}\frac{|\hat{u}_0(\xi)|}{C}\le 0,\any \xi$$
but this can be easily verified after choosing appropriate $C$.\\
\indent The second step is to check that \eqref{inte} decays with respect to $\lam$ for $\lam\ge 1$, but this is quite straightforward. So far the proof of the corollary is finished.$\square$

\begin{rem} 
Instead of \eqref{E:6.1}, what seems to be natual is a scaling like:
\begineq
\v v_0\v_{\hat{L}^{\Phi}}\lesim (\ln \lam)^{\gamma} \v u_0\v_{\hat{L}^{\Phi}}.
\endeq
In order to see that the $\hat{L}^1$ norm is really necessary on the right hand side, we give the following example (for simplicity, take $\gamma=1$ in the definition of $\hat{L}^{\Phi}$),
\begin{displaymath}
\hat{u}_0 (\xi)=\left\{ \begin{array}{ll}								
N& {0\le \xi\le e^{-N}}\\
0 & {\text{else}}
\end{array} \right.
\end{displaymath}
with $N\in \mb{N}$. After some easy calculation, we see that $\v u_0\v_{\hat{L}^{\Phi}}\sim 1$. Then take $\lam= e^{e^{M}}$ with $M\gg N$, define $v_0(x)=\lam u_0(\lam x)$, again it's easy to check that 
\begineq
\v v_0\v_{\hat{L}^{\Phi}}\approx N e^{M}\approx \v \hat{u}_0\v_{L^{\infty}} \ln \lam
\endeq
from which we see how $\v \hat{u}_0\v_{L^{\infty}}$ appears.
\end{rem}

Shaoming Guo, Institute of Mathematics, University of Bonn\\
\indent Address: Endenicher Allee 60, 53115, Bonn\\
\indent Email: shaoming@math.uni-bonn.de

\end{document}